\documentclass[twoside,11pt]{article}
\usepackage[pdftex]{hyperref}
\hypersetup{
 	colorlinks = true,
		allcolors = {red},
}
\usepackage{color}
\usepackage{amsmath}
\usepackage{amssymb}
\usepackage{amsthm}
\usepackage{xparse}
\usepackage{cite}

\newcommand\sa{\smallskipamount}
\newcommand\sLP{\\[\sa]}
\newcommand\sPP{\\[\sa]\indent}
\newcommand\ba{\bigskipamount}
\newcommand\bLP{\\[\ba]}
\newcommand\CC{\mathbb{C}}
\newcommand\RR{\mathbb{R}}
\newcommand\ZZ{\mathbb{Z}}
\newcommand\al\alpha
\newcommand\be\beta
\newcommand\ga\gamma
\newcommand\de\delta
\newcommand\tha\theta
\newcommand\la\lambda
\newcommand\om\omega
\newcommand\Ga{\Gamma}
\newcommand\half{\frac12}
\newcommand\thalf{\tfrac12}
\newcommand\iy\infty
\newcommand\td{\tilde}
\newcommand\wt{\widetilde}
\newcommand\const{{\rm const.}\,}
\newcommand\Zpos{\ZZ_{>0}}
\newcommand\Znonneg{\ZZ_{\ge0}}
\newcommand{\hyp}[5]{\,\mbox{}_{#1}F_{#2}\!\left(
  \genfrac{}{}{0pt}{}{#3}{#4};#5\right)}
\newcommand{\qhyp}[5]{\,\mbox{}_{#1}\phi_{#2}\!\left(
  \genfrac{}{}{0pt}{}{#3}{#4};#5\right)}
\newcommand\LHS{left-hand side}
\newcommand\RHS{right-hand side}
\renewcommand\Re{{\rm Re}\,}
\renewcommand\Im{{\rm Im}\,}
\newcommand{\dup}{\textup{d}}
\newcommand{\eup}{\textup{e}}
\newcommand{\iup}{\mkern1mu\textup{i}\mkern1mu}

\NewDocumentCommand\mycite{m g}{%
  \IfNoValueTF{#2}
    {[\hyperlink{#1}{#1}]}
    {[\hyperlink{#1}{#1}, #2]}%
}
\newcommand\mybibitem[1]{\bibitem[#1]{#1}\hypertarget{#1}{}}

\NewDocumentCommand{\myciteKLS}{m g}{%
  \IfNoValueTF{#2}
    {[\hyperlink{KLS#1}{#1}]}
    {[\hyperlink{KLS#1}{#1}, #2]}%
}
\newcommand\mybibitemKLS[1]{\bibitem[#1]{#1}\hypertarget{KLS#1}{}}

\begin{document}

\title{Additions to the formula lists in
``Hypergeometric orthogonal polynomials and their $q$-analogues''
by Koekoek, Lesky and Swarttouw}
\author{Tom H. Koornwinder}
\date{}
\maketitle
\begin{abstract}
This paper gives a rather arbitrary choice of formulas for
($q$-)hypergeometric orthogonal polynomials belonging to the
($q$-)Askey scheme which the author missed
while consulting Chapters 9 and 14 in the book
``Hypergeometric orthogonal polynomials and their $q$-analogues''
by Koekoek, Lesky and Swarttouw.
The systematics of these chapters will be followed
here, in particular for the numbering of subsections and of references.
\end{abstract}
\subsection*{Introduction}
\label{sec_intro}
This paper contains some formulas for ($q$-)hypergeometric
orthogonal polynomials belonging to the $q$-Askey scheme
which I missed but wanted to use
while consulting Chapters 9 and 14 in the book \mycite{KLS}:
\sLP
R. Koekoek, P.~A. Lesky and R.~F. Swarttouw,
{\em Hypergeometric orthogonal polynomials and their $q$-analogues},
Springer-Verlag, 2010.
\sLP
These chapters form together the (slightly extended) successor of the report
\sLP
R.~Koekoek and  R.~F. Swarttouw,
{\em The Askey-scheme of hypergeometric orthogonal
polynomials and its $q$-analogue},
Report 98-17, Faculty of Technical Mathematics and Informatics,
Delft University of Technology, 1998;
\url{http://aw.twi.tudelft.nl/~koekoek/askey/}.
\sPP
Certainly these chapters give complete lists of formulas of special type,
for instance orthogonality relations and three-term recurrence relations.
But outside these narrow categories there are many other
formulas for ($q$-)orthogonal polynomials which one wants to have available.
Often one can find the desired formula in one of the 
\hyperref[sec_ref1]{standard references} listed at the end of this paper.
Sometimes it is only available in a journal or a less common monograph.
Just for my own comfort, I have brought together some of these formulas.
This will possibly also be helpful for some other users.

Morally, any formula for, say, Jacobi polynomials should fit into a
($q$-)Askey scheme of formulas of such type, where, as always, the arrows
mean specialization or taking limits. However, I have not attempted to reach
such completeness here.
The resulting choice of formulas is
rather arbitrary, just depending on the formulas which I happened to need or
which raised my interest.
For each formula I give  a suitable reference or I sketch a
proof.
It is my intention to gradually extend this collection of formulas.

It should be observed that a more conceptual approach to the ($q$-)Askey scheme
was given smewhat sketchy by Vinet \& Zhedanov \cite{K40} and later
more thoroughly by Verde-Star \cite{K32}. Verde-Star's approach was made
explicit by the author in \cite{K41} for the $q$-Askey scheme,
in \cite{K42} for the Askey scheme, and in \cite{K43} for the Zhedanov type
algebras associated with the $q$-Askey scheme.
Approaches to a $q=-1$ Askey scheme have been made by Pelletier e.a.
\cite{K38} and Verde-Star \cite{K39}. However, orthogonal polynomials
belonging to such a scheme were not concidered in \mycite{KLS} and will
neither be considered here.
\subsection*{Conventions}
\label{sec_conv}
The (x.y) and (x.y.z) type subsection numbers, the
(x.y.z) type formula numbers, and the [x] type citation numbers
refer to \mycite{KLS}.
The (x) type formula numbers refer to this manuscript and the [Kx] type
citation numbers refer to citations which are not in \mycite{KLS}.
Some standard references like \mycite{DLMF}
are given by special acronyms.

$N$ is always a positive integer. Always assume $n$ to be a nonnegative
integer or, if $N$ is present, to be in $\{0,1,\ldots,N\}$.
Throughout assume $0<q<1$.

For each family the coefficient of the term of highest degree of the
orthogonal polynomial of degree $n$ can be found in \mycite{KLS} as the
coefficient of $p_n(x)$ in the formula after the main formula under
the heading ``Normalized Recurrence Relation". If that main formula is numbered
as (x.y.z) then I will refer to the second formula as (x.y.zb).

In the notation of $q$-hypergeometric orthogonal polynomials we
will follow the convention that the parameter list and $q$ are separated
by `$\,|\,$' in the case of a $q$-quadratic lattice (for instance
\hyperref[sec14.1]{Askey--Wilson})
and by `;' in the case of a $q$-linear lattice (for instance
\hyperref[sec14.5]{big $q$-Jacobi}). This convention is mostly followed
in \mycite{KLS}, but not everywhere, see for instance
\hyperref[sec14.20]{little $q$-Laguerre / Wall}.
\subsection*{Acknowledgement}
Many thanks to Howard Cohl for having called my attention so often to
typos and inconsistencies. Thanks also to Roberto Costas Santos
for observing an error and to Gregory Natanson for historical information
about Routh.
\newpage
\subsection*{Contents}
\hyperref[sec_intro]{Introduction}\\
\hyperref[sec_conv]{Conventions}\\
\hyperref[sec_general]{Generalities}
\sLP
\hyperref[sec9.1]{9.1 Wilson}\\
\hyperref[sec9.2]{9.2 Racah}\\
\hyperref[sec9.3]{9.3 Continuous dual Hahn}\\
\hyperref[sec9.4]{9.4 Continuous Hahn}\\
\hyperref[sec9.5]{9.5 Hahn}\\
\hyperref[sec9.6]{9.6 Dual Hahn}\\
\hyperref[sec9.7]{9.7 Meixner--Pollaczek}\\
\hyperref[sec9.8]{9.8 Jacobi}

\hyperref[sec9.8.1]{9.8.1 Gegenbauer / Ultraspherical}

\hyperref[sec9.8.2]{9.8.2 Chebyshev}\\
\hyperref[sec9.9]{9.9 Pseudo Jacobi (or Romanovski--Routh)}\\
\hyperref[sec9.10]{9.10 Meixner}\\
\hyperref[sec9.11]{9.11 Krawtchouk}\\
\hyperref[sec9.12]{9.12 Laguerre}\\
\hyperref[sec9.14]{9.13 Bessel}\\
\hyperref[sec9.14]{9.14 Charlier}\\
\hyperref[sec9.15]{9.15 Hermite}
\sLP
\hyperref[sec14.1]{14.1 Askey--Wilson}\\
\hyperref[sec14.2]{14.2 $q$-Racah}\\
\hyperref[sec14.3]{14.3 Continuous dual $q$-Hahn}\\
\hyperref[sec14.4]{14.4 Continuous $q$-Hahn}\\
\hyperref[sec14.5]{14.5 Big $q$-Jacobi}\\
\hyperref[sec14.7]{14.7 Dual $q$-Hahn}\\
\hyperref[sec14.8]{14.8 Al-Salam--Chihara}\\
\hyperref[sec14.9]{14.9 $q$-Meixner--Pollaczek}\\
\hyperref[sec14.10]{14.10 Continuous $q$-Jacobi}

\hyperref[sec14.10.1]{14.10.1 Continuous $q$-ultraspherical / Rogers}\\
\hyperref[sec14.11]{14.11 Big $q$-Laguerre}\\
\hyperref[sec14.12]{14.12 Little $q$-Jacobi}\\
\hyperref[sec14.14]{14.14 Quantum $q$-Krawtchouk}\\
\hyperref[sec14.16]{14.16 Affine $q$-Krawtchouk}\\
\hyperref[sec14.17]{14.17 Dual $q$-Krawtchouk}\\
\hyperref[sec14.20]{14.20 Little $q$-Laguerre / Wall}\\
\hyperref[sec14.21]{14.21 $q$-Laguerre}\\
\hyperref[sec14.27]{14.27 Stieltjes-Wigert}\\
\hyperref[sec14.28]{14.28 Discrete $q$-Hermite I}\\
\hyperref[sec14.29]{14.29 Discrete $q$-Hermite II}
\sLP
\hyperref[sec_ref1]{Standard references}\\
\hyperref[sec_ref2]{References from [KLS]}\\
\hyperref[sec_ref3]{Other references}
\newpage
\subsection*{Generalities}
\label{sec_general}
\paragraph{Criteria for uniqueness of orthogonality measure}\;\\
According to Shohat \& Tamarkin \cite[p.50]{K6}
orthonormal polynomials $p_n$ have a unique orthogonality measure (up to positive
constant factor) if
for some $z\in\CC$ we have
\begin{equation}
\sum_{n=0}^\iy |p_n(z)|^2 = \iy.
\label{90}
\end{equation}

Also (see Shohat \& Tamarkin \cite[p.59]{K6}),
monic orthogonal polynomials $p_n$ with three-term recurrence relation
$x p_n(x) = p_{n+1}(x)+B_n p_n(x)+C_n p_{n-1}(x)$
($C_n$ necessarily positive)
have a unique orthogonality measure if
\begin{equation}
\sum_{n=1}^\iy (C_n)^{-1/2}=\iy.
\label{93}
\end{equation}

Furthermore, if orthogonal polynomials have an orthogonality measure with
bounded support, then this is unique (see Chihara \myciteKLS{146}).
\paragraph{Kernel polynomials and the three-term recurrence relation}\;\\
For given monic orthogonal polynomials $\{p_n\}$ with respect to orthogonality
measure $\mu$ and with
\begin{equation*}
h_n:=\int_\RR p_n(x)^2\,\dup\mu(x),
\end{equation*}
there is the \emph{Christoffel--Darboux formula}
\begin{equation}
K_n(x,y):=\sum_{k=0}^n \frac{p_k(x)p_k(y)}{h_k}
=\frac1{h_n}\,\frac{p_{n+1}(x)p_n(y)-p_n(x)p_{n+1}(y)}{x-y}\qquad(x\ne y).
\label{205}
\end{equation}
Fix $y\in\RR$ and suppose that $\operatorname{supp}(\mu)\subseteq(-\iy,y]$.
Then $p_n(y)\ne0$ for all $n$ and the monic polynomials
\begin{equation}
q_n(x):=\frac{h_n}{p_n(y)}\,K_n(x,y)
\label{206}
\end{equation}
are orthogonal with respect to $(y-x)\,\dup\mu(x)$. They are called
\emph{kernel polynomials} (see Chihara \myciteKLS{146}{Ch.~1, \S7}).
There is a pair of contiguous relations relating the polynomialsd $p_n$
and $q_n$:
\begin{align}
(x-y)q_n(x)&=p_{n+1}(x)-A_n p_n(x),\label{207}\\
p_n(x)&=q_n(x)-C_n q_{n-1}(x),\label{208}
\end{align}
where
\begin{equation}
A_n=\frac{p_{n+1}(y)}{p_n(y)}\,,\qquad
C_n=\frac{h_n}{h_{n-1}}\,\frac{p_{n-1}(y)}{p_n(y)}\,.
\label{209}
\end{equation}
Then the \emph{three-term recurrence relations} for the orthogonal polynomials
$p_n$ and $q_n$ can be written in the form
(see \cite[\S5, Lemma 1]{K33})
\begin{align}
x\,p_n(x)&=p_{n+1}(x)+(y-A_n-C_n)p_n(x)+A_{n-1}C_n p_{n-1}(x),\label{210}\\
x\,q_n(x)&=q_{n+1}(x)+(y-A_n-C_{n+1})q_n(x)+A_nC_n q_{n-1}(x).\label{211}
\end{align}
In the above formulas put terms containing the factor $C_0$ equal to 0.

In many cases in \mycite{KLS}{Chapters 9, 14} the normalized three-term
recurrence relation is given in the form \eqref{210}, already in the
Askey--Wilson case (14.1.5), and where it is not written in this way, it can be
done so. See for instance \eqref{200} for Jacobi.

If we write the normalized recurrence relation for the $p_n$ as
\begin{equation}
x\,p_n(x)=p_{n+1}(x)+b_n\,p_n(x)+c_n\,p_{n-1}(x),
\label{212}
\end{equation}
and compare it with \eqref{210} then
\begin{equation}
b_0=y-A_0,\qquad
b_n=y-A_n-C_n,\qquad
c_n=A_{n-1}C_n\qquad (n\ge 1).
\label{213}
\end{equation}
This can be recursively solved for the $A_n,C_n$ in terms of the $b_n,c_n$ by
\begin{equation}
A_0=y-b_0,\qquad
C_n=\frac{c_n}{A_{n-1}}\,,\quad
A_n=y-b_n-C_n\qquad (n\ge1).
\label{214}
\end{equation}

Equations \eqref{207}, \eqref{208}, \eqref{210} correspond to an LU
factorization of the Jacobi matrix associated with the OPs $p_n$, see
\cite[Lemma 2.1]{K34}, where also \eqref{214} is given.
\paragraph{Even orthogonality measure}\;\\
If $\{p_n\}$ is a system of orthogonal polynomials with respect to an even
orthogonality measure which satisfies the three-term recurrence relation
\begin{equation*}
x p_n(x)=a_n p_{n+1}(x)+c_n\,p_{n-1}(x)
\end{equation*}
then
\begin{equation}
\frac{p_{2n}(0)}{p_{2n-2}(0)}=-\,\frac{c_{2n-1}}{a_{2n-1}}\,.
\label{1}
\end{equation}
\paragraph{Finite systems of OPs of degree up to $N$ with weights on $N+1$
points}\;\\
Suppose we have OPs $\{p_n\}_{n=0}^N$ which are orthogonal on
$\{x_0,x_1,\ldots,x_N\}$ with respect to weights $w_i$ ($i=0,1,\ldots,N$).
Then we have recurrence relations
\begin{equation}
xp_n(x)=A_n p_{n+1}(x)+B_n p_n(x)+C_n p_{n-1}(x)\quad(n=0,1\ldots,N),
\label{221}
\end{equation}
where $p_{-1}(x)=0$, $p_{N+1}(x)=(x-x_0)\ldots(x-x_N)$ and $p_N(x)=A_N x^N+
\mbox{terms of lower degree}$.
For a proof of the case $n=N$ note that, for $x\in\{x_0,x_1,\ldots,x_N\}$,
we have for certain coefficients $B_N,C_N$ that
$xp_N(x)=B_N p_N(x)+C_N p_{N-1}(x)$ by orthogonality and by the fact that
$p_0,p_1,\ldots,p_N$ is a basis of the function space on this set. Hence
$xp_N(x)-B_N p_N(x)-C_N p_{M-1}(x)$ is a polynomial of degree $N+1$
which vanishes
on $\{x_0,x_1,\ldots,x_N\}$ and for which the coefficient of $x^{N+1}$ equals
the coefficient of $x^N$ for $p_N(x)$.
Hence $xp_N(x)-B_N p_N(x)-C_N p_{N-1}(x)=A_N(x-x_0)\ldots(x-x_N)$.
\paragraph{Zeros of an OP of degree $n$}\label{generalities_zeros}\;\\
Let $\mu$ be a positive measure with infinite support on a closed interval $I$
such that $|x|^{2n-1}$ (but not necessarily $|x|^{2n}$) is integrable on $I$
with respect to $\mu$. Let $p_n(x)$ be a polynomial of degree~$n$ such that
$\int_I p_n(x)\,x^k\,\dup\mu(x)=0$ for $k=0,1,\ldots,n-1$. Then $p_n(x)$
has $n$ distinct zeros on the interior of $I$.
\sLP
\textbf{Proof}
Suppose that $p_n$ has precisely $k<n$ sign changes on points $x_1,\ldots,x_k$
in the interior of~$I$. Then
$0=\int_I p_n(x)\,\prod_{j=1}^k(x-x_j)\,\dup\mu(x)\ne0$, on the one hand by
orthogonality and on the other hand because the integrand has a definite
sign outside $x_1,\ldots,x_k$. So we have a contradiction.\qed
\paragraph{Appell's bivariate hypergeometric function $F_4$}
This is defined by
\begin{equation}
F_4(a,b;c,c';x,y):=\sum_{m,n=0}^\iy\frac{(a)_{m+n}(b)_{m+n}}{(c)_m(c')_n\,m!\,n!}\,
x^my^n\qquad(|x|^\half+|y|^\half<1),
\label{62}
\end{equation}
see \mycite{HTF1}{5.7(9), 5.7(44)} or \mycite{DLMF}{(16.13.4)}.
There is the reduction formula
\begin{equation*}
F_4\left(a,b;b,b;\frac{-x}{(1-x)(1-y)},\frac{-y}{(1-x)(1-y)}\right)=
(1-x)^a(1-y)^a\,\hyp21{a,1+a-b}b{xy},
\end{equation*}
see \mycite{HTF1}{5.10(7)}. When combined with the quadratic transformation
\mycite{HTF1}{2.11(34)} (here $a-b-1$ should be replaced by $a-b+1$),
see also \mycite{DLMF}{(15.8.15)}, this yields
\begin{multline*}
F_4\left(a,b;b,b;\frac{-x}{(1-x)(1-y)},\frac{-y}{(1-x)(1-y)}\right)\\
=\left(\frac{(1-x)(1-y)}{1+xy}\right)^a\,
\hyp21{\thalf a,\thalf(a+1)}b{\frac{4xy}{(1+xy)^2}}.
\end{multline*}
This can be rewritten as
\begin{equation}
F_4(a,b;b,b;x,y)=(1-x-y)^{-a}\,\hyp21{\thalf a,\thalf(a+1)}b
{\frac{4xy}{(1-x-y)^2}}.
\label{63}
\end{equation}
Note that, if $x,y\ge0$ and $x^\half+y^\half<1$, then
$1-x-y>0$ and $0\le\frac{4xy}{(1-x-y)^2}<1$.
\paragraph{$q$-Hypergeometric series of base $q^{-1}$}
By \mycite{GR}{Exercise 1.4(i)}:
\begin{equation}
\qhyp rs{a_1,\ldots,a_r}{b_1,\ldots b_s}{q^{-1},z}
=\qhyp{s+1}s{a_1^{-1},\ldots a_r^{-1},0,\ldots,0}
{b_1^{-1},\ldots,b_s^{-1}}{q,\frac{qa_1\ldots a_rz}{b_1\ldots b_s}}
\label{154}
\end{equation}
for $r\le s+1$, $a_1,\ldots,a_r,b_1,\ldots,b_s\ne0$.
In the non-terminating case, for $0<q<1$, there is convergence if
$|z|<b_1\ldots b_s/(qa_1\ldots a_r)$\,.
\paragraph{A transformation of a terminating ${}_2\phi_1$}
By \mycite{GR}{Exercise 1.15(i)} we have
\begin{equation}
\qhyp21{q^{-n},b}c{q,z}=(bz/(cq);q^{-1})_n\,
\qhyp32{q^{-n},c/b,0}{c,cq/(bz)}{q,q}.
\label{151}
\end{equation}
\paragraph{Very-well-poised $q$-hypergeometric series}
The notation of \mycite{GR}{(2.1.11)} will be followed:
\begin{equation}
{}_{r+1}W_r(a_1;a_4,a_5,\ldots,a_{r+1};q,z):=
\qhyp{r+1}r{a_1,qa_1^\half,-qa_1^\half,a_4,\ldots,a_{r+1}}
{a_1^\half,-a_1^\half,qa_1/a_4,\ldots,qa_1/a_{r+1}}{q,z}.
\label{111}
\end{equation}
\paragraph{Theta function}
The notation of \mycite{GR}{(11.2.1)} will be followed:
\begin{equation}
\tha(x;q):=(x,q/x;q)_\iy,\qquad
\tha(x_1,\ldots,x_m;q):=\tha(x_1;q)\ldots\tha(x_m;q).
\label{117}
\end{equation}
\subsection*{9.1 Wilson}
\label{sec9.1}
\paragraph{Symmetry}
The Wilson polynomial $W_n(y;a,b,c,d)$ is symmetric
in $a,b,c,d$.
\\
This follows from the orthogonality relation (9.1.2)
together with the value of its coefficient of $y^n$ given in (9.1.5b).
Alternatively, combine (9.1.1) with \mycite{AAR}{Theorem 3.1.1}.\\
As a consequence, it is sufficient to give generating function (9.1.12). Then the generating
functions (9.1.13), (9.1.14) will follow by symmetry in the parameters.
\paragraph{Hypergeometric representation}
In addition to (9.1.1) we have (see \myciteKLS{513}{(2.2)}):
\begin{multline}
W_n(x^2;a,b,c,d)
=\frac{(a-\iup x)_n (b-\iup x)_n (c-\iup x)_n (d-\iup x)_n}{(-2\iup x)_n}\\
\times\hyp76{2\iup x-n,\iup x-\thalf n+1,a+\iup x,b+\iup x,c+\iup x,d+\iup x,-n}
{\iup x-\thalf n,1-n-a+\iup x,1-n-b+\iup x,1-n-c+\iup x,1-n-d+\iup x,1+2\iup x}
1.
\label{112}
\end{multline}
The symmetry in $a,b,c,d$ is clear from \eqref{112}.
\paragraph{Special value}
\begin{equation}
W_n(-a^2;a,b,c,d)=(a+b)_n(a+c)_n(a+d)_n\,,
\label{91}
\end{equation}
and similarly for arguments $-b^2$, $-c^2$ and
$-d^2$ by symmetry of $W_n$ in $a,b,c,d$.
\paragraph{Uniqueness of orthogonality measure}
Under the assumptions on $a,b,c,d$ for (9.1.2) or (9.1.3) the orthogonality
measure is unique up to constant factor.

For the proof assume without
loss of generality (by the symmetry in $a,b,c,d$) that $\Re a\ge0$.
Write the \RHS\ of (9.1.2) or (9.1.3) as $h_n\de_{m,n}$.
Observe from (9.1.2) and \eqref{91} that
\[
\frac{|W_n(-a^2;a,b,c,d)|^2}{h_n} = O(n^{4\Re a-1})\quad\hbox{as $n\to\iy$.}
\]
Therefore \eqref{90} holds, from which the uniqueness of the orthogonality
measure follows.

By a similar, but necessarily more complicated argument Ismail et al.\
\myciteKLS{281}{Section 3} proved the uniqueness of orthogonality measure for
associated Wilson polynomials.
\subsection*{9.2 Racah}
\label{sec9.2}
\paragraph{Racah in terms of Wilson}
In the Remark on p.196 Racah polynomials are expressed in terms of
Wilson polynomials. This can be equivalently written as
\begin{multline}
R_n\big(x(x-N+\de);\al,\be,-N-1,\de\big)\\
=\frac{W_n\big(-(x+\thalf(\de-N))^2;\thalf(\de-N),\al+1-\thalf(\de-N),
\be+\thalf(\de+N)+1,-\half(\de+N)\big)}
{(\al+1)_n (\be+\de+1)_n (-N)_n}\,.
\label{146}
\end{multline}
\subsection*{9.3 Continuous dual Hahn}
\label{sec9.3}
\paragraph{Symmetry}
The continuous dual Hahn polynomial $S_n(y;a,b,c)$ is symmetric
in $a,b,c$.\\
This follows from the orthogonality relation (9.3.2)
together with the value of its coefficient of $y^n$ given in (9.3.5b).
Alternatively, combine (9.3.1) with \mycite{AAR}{Corollary 3.3.5}.\\
As a consequence, it is sufficient to give generating function (9.3.12). Then the generating
functions (9.3.13), (9.3.14) will follow by symmetry in the parameters.
\paragraph{Special value}
\begin{equation}
S_n(-a^2;a,b,c)=(a+b)_n(a+c)_n\,,
\label{92}
\end{equation}
and similarly for arguments $-b^2$ and $-c^2$ by symmetry of $S_n$ in $a,b,c$.
\paragraph{Uniqueness of orthogonality measure}
Under the assumptions on $a,b,c$ for (9.3.2) or (9.3.3) the orthogonality
measure is unique up to constant factor.

For the proof assume without
loss of generality (by the symmetry in $a,b,c$) that $\Re a\geq0$.
Write the \RHS\ of (9.3.2) or (9.3.3) as $h_n\de_{m,n}$.
Observe from (9.3.2) and \eqref{92} that
\[
\frac{|S_n(-a^2;a,b,c)|^2}{h_n} = O(n^{2\Re a-1})\quad
\hbox{as $n\to\iy$.}
\]
Therefore \eqref{90} holds, from which the uniqueness of the orthogonality
measure follows.
\paragraph{Special continuous dual Hahn in terms of Wilson}
\begin{equation}
S_n\big(x;a,b,\thalf\big)=\frac{2^{2n}}{(a+b+n)_n}\,
W_n\big(\tfrac14 x;\thalf a,\thalf(a+1),\thalf b,\thalf(b+1)\big).
\label{192}
\end{equation}
For the proof compare the weight functions and the values for $x=-a^2$.
\paragraph{Generating functions}
By (9.3.17) the generating function (9.3.16) has the generating function
(9.7.13) for Meixner--Pollaczek polynomials as a limit case.
\subsection*{9.4 Continuous Hahn}
\label{sec9.4}
\paragraph{Orthogonality relation and parameter symmetry}
The orthogonality relation (9.4.2) holds under the more general assumption that
$\Re(a,b,c,d)>0$ and $(c,d)=(\overline a,\overline b)$ or $(\overline b,\overline a)$.\\
Thus, under these assumptions, the continuous Hahn polynomial
$p_n(x;a,b,c,d)$
is symmetric in $a,b$ and in $c,d$.
This follows from the orthogonality relation (9.4.2)
together with the value of its coefficient of $x^n$ given in (9.4.4b).\\
As a consequence, it is sufficient to give generating function (9.4.11). Then the generating
function (9.4.12) will follow by symmetry in the parameters.
\paragraph{Symmetry}
\begin{equation}
p_n(-x;a,b,\overline a,\overline b)
=(-1)^n p_n(x;\overline a,\overline b,a,b).
\end{equation}
\paragraph{Special value}
\begin{equation}
p_n(\iup a;a,b,\overline a,\overline b)=
\frac{\iup^n(a+\overline a)_n (a+\overline b)_n}{n!}\,.
\end{equation}
Similarly, $p_n(x;a,b,\overline a,\overline b)$ has special values
for $x=-\iup \overline a$, $\iup b$ and $-\iup \overline b$.
\paragraph{Quadratic transformation}
For $a,b\in\RR$ or $b=\overline a$ we have \cite[(2.29), (2.30)]{K30}
\begin{equation}
\frac{p_{2n}(x;a,b,\overline a,\overline b)}
{p_{2n}(\iup a;a,b,\overline a,\overline b)}=
\frac{W_n(x^2;a,b,\thalf,0)}{W_n(-a^2;a,b,\thalf,0)}\,,\quad
\frac{p_{2n+1}(x;a,b,\overline a,\overline b)}
{p_{2n+1}(\iup a;a,b,\overline a,\overline b)}=
\frac{x W_n(x^2;a,b,\thalf,1)}{\iup a W_n(-a^2;a,b,\thalf,1)}\,.
\label{187}
\end{equation}
\paragraph{Explicit expression}
For $a,b\in\RR$ or $b=\overline a$ we have by \eqref{187}, (9.1.1)
and reversion of direction of summation that
\begin{multline}
p_n(x;a,b,\overline a,\overline b)
=\frac{(n+a+b+\overline a+\overline b-1)_n}{n!}\,
x^{n-2[\half n]}\,(-\thalf n+\iup x+1)_{[\half n]}
(-\thalf n-\iup x+1)_{[\half n]}\\
\times\hyp43{-\thalf n,-\thalf n+\thalf,-\thalf n-a+1,-\thalf n-b+1}
{-n-a-b+\tfrac32,-\thalf n+\iup x+1,-\thalf n-\iup x+1}1.
\end{multline}
\paragraph{Special cases}
In the following special case there is a reduction to
Meixner--Pollaczek:
\begin{equation}
p_n(x;a,a+\thalf,a,a+\thalf)=
\frac{(2a)_n (2a+\thalf)_n}{(4a)_n}\,P_n^{(2a)}(2x;\thalf\pi).
\end{equation}
See \myciteKLS{342}{(2.6)} (note that in \myciteKLS{342}{(2.3)} the
Meixner--Pollaczek polyonmials are defined different from (9.7.1),
without a constant factor in front).

For $0<a<1$ the continuous Hahn polynomials $p_n(x;a,1-a,a,1-a)$
are orthogonal on $(-\iy,\iy)$ with respect to the weight function
$\big(\cosh(2\pi x)-\cos(2\pi a)\big)^{-1}$
(by straightforward computation from (9.4.2)).
For $a=\tfrac14$ the two special cases coincide:
Meixner--Pollaczek with weight function $\big(\cosh(2\pi x)\big)^{-1}$.
\paragraph{Uniqueness of orthogonality measure}
The coefficient of $p_{n-1}(x)$ in (9.4.4) behaves as $O(n^2)$ as $n\to\iy$.
Hence \eqref{93} holds, by which the orthogonality measure is unique.
\subsection*{9.5 Hahn}
\label{sec9.5}
\paragraph{Special values}
\begin{equation}
Q_n(0;\al,\be,N)=1,\quad
Q_n(N;\al,\be,N)=\frac{(-1)^n(\be+1)_n}{(\al+1)_n}\,.
\label{95}
\end{equation}
Use (9.5.1) and compare with (9.8.1) and \eqref{50}.

From (9.5.3) and \eqref{1} it follows that
\begin{equation}
Q_{2n}(N;\al,\al,2N)=\frac{(\thalf)_n(N+\al+1)_n}{(-N+\thalf)_n(\al+1)_n}\,.
\label{30}
\end{equation}
From (9.5.1) and \mycite{DLMF}{(15.4.24)} it follows that
\begin{equation}
Q_N(x;\al,\be,N)=\frac{(-N-\be)_x}{(\al+1)_x}\qquad(x=0,1,\ldots,N).
\label{44}
\end{equation}
\paragraph{Symmetries}
By the orthogonality relation (9.5.2):
\begin{equation}
\frac{Q_n(N-x;\al,\be,N)}{Q_n(N;\al,\be,N)}=Q_n(x;\be,\al,N),
\label{96}
\end{equation}
It follows from \eqref{97} and \eqref{45} that
\begin{equation}
\frac{Q_{N-n}(x;\al,\be,N)}{Q_N(x;\al,\be,N)}
=Q_n(x;-N-\be-1,-N-\al-1,N)
\qquad(x=0,1,\ldots,N).
\label{100}
\end{equation}
\paragraph{Duality}
The Remark on p.208 gives the duality between Hahn and dual Hahn polynomials:
\begin{equation}
Q_n(x;\al,\be,N)=R_x(n(n+\al+\be+1);\al,\be,N)\quad(n,x\in\{0,1,\ldots N\}).
\label{45}
\end{equation}
\subsection*{9.6 Dual Hahn}
\label{sec9.6}
\paragraph{Special values}
By \eqref{44} and \eqref{45} we have
\begin{equation}
R_n(N(N+\ga+\de+1);\ga,\de,N)=\frac{(-N-\de)_n}{(\ga+1)_n}\,.
\label{47}
\end{equation}
It follows from \eqref{95} and \eqref{45} that
\begin{equation}
R_N(x(x+\ga+\de+1);\ga,\de,N)
=\frac{(-1)^x(\de+1)_x}{(\ga+1)_x}\qquad(x=0,1,\ldots,N).
\label{101}
\end{equation}
\paragraph{Symmetries}
Write the weight in (9.6.2) as
\begin{equation}
w_x(\al,\be,N):=N!\,\frac{2x+\ga+\de+1}{(x+\ga+\de+1)_{N+1}}\,
\frac{(\ga+1)_x}{(\de+1)_x}\,\binom Nx.
\label{98}
\end{equation}
Then
\begin{equation}
(\de+1)_N\,w_{N-x}(\ga,\de,N)=
(-\ga-N)_N\,w_x(-\de-N-1,-\ga-N-1,N).
\label{99}
\end{equation}
Hence, by (9.6.2),
\begin{equation}
\frac{R_n((N-x)(N-x+\ga+\de+1);\ga,\de,N)}{R_n(N(N+\ga+\de+1);\ga,\de,N)}
=R_n(x(x-2N-\ga-\de-1);-N-\de-1,-N-\ga-1,N).
\label{97}
\end{equation}
Alternatively, \eqref{97} follows from (9.6.1) and
\mycite{DLMF}{(16.4.11)}.

It follows from \eqref{96} and \eqref{45} that
\begin{equation}
\frac{R_{N-n}(x(x+\ga+\de+1);\ga,\de,N)}
{R_N(x(x+\ga+\de+1);\ga,\de,N)}
=R_n(x(x+\ga+\de+1);\de,\ga,N)\qquad(x=0,1,\ldots,N).
\label{102}
\end{equation}
\paragraph{Re: (9.6.11).}
The generating function (9.6.11) can be written in a more conceptual way as
\begin{equation}
(1-t)^x\,\hyp21{x-N,x+\ga+1}{-\de-N}t=\frac{N!}{(\de+1)_N}\,
\sum_{n=0}^N \om_n\,R_n(\la(x);\ga,\de,N)\,t^n,
\label{2}
\end{equation}
where
\begin{equation}
\om_n:=\binom{\ga+n}n \binom{\de+N-n}{N-n},
\label{3}
\end{equation}
i.e., the denominator on the \RHS\ of (9.6.2).
By the duality between Hahn polynomials and dual Hahn polynomials (see \eqref{45}) the above generating function can be rewritten in
terms of Hahn polynomials:
\begin{equation}
(1-t)^n\,\hyp21{n-N,n+\al+1}{-\be-N}t=\frac{N!}{(\be+1)_N}\,
\sum_{x=0}^N w_x\,Q_n(x;\al,\be,N)\,t^x,
\label{4}
\end{equation}
where
\begin{equation}
w_x:=\binom{\al+x}x \binom{\be+N-x}{N-x},
\label{5}
\end{equation}
i.e., the weight occurring in the orthogonality relation (9.5.2)
for Hahn polynomials.
\paragraph{Re: (9.6.15).}
There should be a closing bracket before the equality sign.
\subsection*{9.7 Meixner--Pollaczek}
\label{sec9.7}
\paragraph{Re: (9.7.1)}
In addition to the hypergeometric representation (9.7.1) we have, by
the Pfaff transformation \mycite{HTF1}{2.9(3)}, that
\begin{equation}
P_n^{(\la)}(x;\phi)=\frac{(2\la)_n}{n!}\,\eup^{-\iup n\phi}
\hyp21{-n,\la-\iup x}{2\la}{1-\eup^{2\iup\phi}}.
\label{195}
\end{equation}
\paragraph{Special values}
By (9.7.1) and \eqref{195} we have:
\begin{equation}
P_n^{(\la)}(\iup\la;\phi)=\frac{(2\la)_n}{n!}\,\eup^{\iup n\phi},\qquad
P_n^{(\la)}(-\iup\la;\phi)=\frac{(2\la)_n}{n!}\,\eup^{-\iup n\phi}.
\label{193}
\end{equation}
\paragraph{Symmetry}
\begin{equation}
P_n^{(\la)}(x;\phi)=(-1)^n P_n^{(\la)}(-x;\pi-\phi).
\label{196}
\end{equation}

\paragraph{Quadratic transformations}
\cite[(2.33), (2.34)]{K30}
\begin{equation}
\frac{P_{2n}^{(a)}(x;\thalf\pi)}{P_{2n}^{(a)}(\iup a;\thalf\pi)}=
\frac{S_n(x^2;a,\thalf,0)}{S_n(-a^2;a,\thalf,0)}\,,\qquad
\frac{P_{2n+1}^{(a)}(x;\thalf\pi)}{P_{2n+1}^{(a)}(\iup a;\thalf\pi)}=
\frac{x S_n(x^2;a,\thalf,1)}{\iup aS_n(-a^2;a,\thalf,1)}\,.\label{194}
\end{equation}
These are limit cases of \eqref{187} by the limits (9.1.16), (9.4.14).
\paragraph{Uniqueness of orthogonality measure}
The coefficient of $p_{n-1}(x)$ in (9.7.4) behaves as $O(n^2)$ as $n\to\iy$.
Hence \eqref{93} holds, by which the orthogonality measure is unique.
\paragraph{Generating functions}
By (9.3.17) the generating function (9.3.16) for continuous dual Hahn
polynomials has the generating function
(9.7.13) as a limit case. By (9.7.14) formula (9.7.13) has the generating
function (9.12.12) for Laguerre polynomials as a limit case.
\subsection*{9.8 Jacobi}
\label{sec9.8}
\paragraph{Orthogonality relation}
Write the \RHS\ of (9.8.2) as $h_n\,\de_{m,n}$. Then
\begin{equation}
\begin{split}
&\frac{h_n}{h_0}=
\frac{n+\al+\be+1}{2n+\al+\be+1}\,
\frac{(\al+1)_n(\be+1)_n}{(\al+\be+2)_n\,n!}\,,\quad
h_0=\frac{2^{\al+\be+1}\Ga(\al+1)\Ga(\be+1)}{\Ga(\al+\be+2)}\,,\sLP
&\frac{h_n}{h_0\,(P_n^{(\al,\be)}(1))^2}=
\frac{n+\al+\be+1}{2n+\al+\be+1}\,
\frac{(\be+1)_n\,n!}{(\al+1)_n\,(\al+\be+2)_n}\,.
\end{split}
\label{60}
\end{equation}

In (9.8.3) the numerator factor $\Ga(n+\al+\be+1)$ in the last line should be
$\Ga(\be+1)$. When thus corrected, (9.8.3) can be rewritten as:
\begin{equation}
\begin{split}
&\int_1^\iy P_m^{(\al,\be)}(x)\,P_n^{(\al,\be)}(x)\,(x-1)^\al (x+1)^\be\,
\dup x=h_n\,\de_{m,n}\,,\\
&\qquad\qquad\qquad\qquad\qquad\qquad\qquad\quad-1-\be>\al>-1,\quad m,n<-\thalf(\al+\be+1),\\
&\frac{h_n}{h_0}=
\frac{n+\al+\be+1}{2n+\al+\be+1}\,
\frac{(\al+1)_n(\be+1)_n}{(\al+\be+2)_n\,n!}\,,\quad
h_0=\frac{2^{\al+\be+1}\Ga(\al+1)\Ga(-\al-\be-1)}{\Ga(-\be)}\,.
\end{split}
\label{122}
\end{equation}
Following Lesky \myciteKLS{382} the Jacobi polynomials in case of
orthogonality relation \eqref{122} may be called
\emph{Romanovski--Jacobi polynomials}.

The orthogonality \eqref{122} remains true if $-1-\be>\al>-1$ and
$m<n<-\half(\al+\be)$.
Hence for $0<n<-\half(\al+\be)$ the polynomial
${}_2F_1(-n,n+\al+\be+1;\al+1;x)=\const P_n^{(\al,\be}(1-2x)=
\const P_n^{(\be,\al)}(2x-1)$ has all its $n$ zeros
on the interval $(-\iy,0)$ if $\al>-1$, and on $(1,\iy)$ if $\be>-1$,
see p.\pageref{generalities_zeros}.
Hence the polynomial ${}_2F_1(-n,b;c;x)$ fhs all its $n$ zeros on
$(-\iy,0)$ if $b<-n+1$ and $c>0$ and on
$(1,\iy)$ if $b<-n+1$ and $b-c-n+1>0$, see \cite[Theorem 1.1]{K37}.
\paragraph{Symmetry}
\begin{equation}
P_n^{(\al,\be)}(-x)=(-1)^n\,P_n^{(\be,\al)}(x).
\label{48}
\end{equation}
Use (9.8.2) and (9.8.5b) or see \mycite{DLMF}{Table 18.6.1}.
\paragraph{Special values}
\begin{equation}
P_n^{(\al,\be)}(1)=\frac{(\al+1)_n}{n!}\,,\quad
P_n^{(\al,\be)}(-1)=\frac{(-1)^n(\be+1)_n}{n!}\,,\quad
\frac{P_n^{(\al,\be)}(-1)}{P_n^{(\al,\be)}(1)}=\frac{(-1)^n(\be+1)_n}{(\al+1)_n}\,.
\label{50}
\end{equation}
Use (9.8.1) and \eqref{48} or see \mycite{DLMF}{Table 18.6.1}.
\paragraph{Normalized recurrence relation}
Formula (9.8.5) can be rewritten as
\begin{equation}
x\,p_n(x)=p_{n+1}(x)+(1-A_n-C_n)p_n(x)+A_{n-1}C_n\,p_{n-1}(x),
\label{200}
\end{equation}
where $p_n(x)=2^n n!\,P_n^{(\al,\be)}(x)/(n+\al+\be+1)_n$ and
\[
A_n=\frac{2(n+\al+1)(n+\al+\be+1)}{(2n+\al+\be+1)(2n+\al+\be+2)}\,,\qquad
C_n=\frac{2n(n+\be)}{(2n+\al+\be)(2n+\al+\be+1)}\,.
\]
\paragraph{Contiguous relations}
\begin{align}
&(n+\thalf\al+\thalf\be+1)(1-x)P_n^{(\al+1,\be)}(x)
=-(n+1)P_{n+1}^{(\al,\be)}(x)+(n+\al+1)P_n^{(\al,\be)}(x),
\label{201}\\
&(2n+\al+\be+1)P_n^{(\al,\be)}(x)
=(n+\al+\be+1)P_n^{(\al+1,\be)}(x)-(n+\be)P_{n-1}^{(\al+1,\be)}(x).
\label{202}
\end{align}
See \mycite{HTF2}{10.8(32) and (35)}. These can be rewritten as
\begin{align}
(x-1)q_n(x)&=p_{n+1}(x)-A_n p_n(x),\label{203}\\
p_n(x)&=q_n(x)-C_n q_{n-1}(x),\label{204}
\end{align}
where $q_n(x)=2^n n!\,P_n^{(\al+1,\be)}(x)/(n+\al+\be+2)_n$ and
$p_n(x)$, $A_n$ and $C_n$ are as above.

Formula \eqref{200} can be derived from \eqref{203}, \eqref{204} by substituting
these last two formulas in the following rewritten form of \eqref{200}
(compare with \eqref{207}--\eqref{210}):
\[
(x-1)p_n(x)=\big(p_{n+1}(x)-A_n p_n(x)\big)
-C_n\big(p_n(x)-A_{n-1} p_{n-1}(x)\big).
\]
\paragraph{Generating functions}
Formula (9.8.15) was first obtained by Brafman \myciteKLS{109}{(12)}.
Alternatively (see \myciteKLS{109}{(9)} or use
\mycite{DLMF}{(16.16.6)}), the \LHS\ of (9.8.15)
can be written as Appell's hypergeometric function $F_4$:
\begin{equation}
F_4\big(\ga,\al+\be+1-\ga;\al+1,\be+1;\thalf t(x-1),\thalf t(x+1)\big)
=\sum_{k=0}^\iy \frac{(\ga)_k (\al+\be+1-\ga)_k}{(\al+1)_k (\be+1)_k}\,
P_k^{(\al,\be)}(x) t^k
\label{198}
\end{equation}
The generating function (9.12.12) for Laguerre polynomials
is a limit case of \eqref{198} by (9.8.16).

Formula (9.8.15) with $t$, $x$ replaced by $\thalf(x+y)$,
$\frac{1+xy}{x+y}$, respectively, takes the form
\begin{multline}
\hyp21{\ga,\al+\be+1-\ga}{\al+1}{\thalf(1-x)}
\hyp21{\ga,\al+\be+1-\ga}{\be+1}{\thalf(1+y)}\\
=\sum_{k=0}^\iy \frac{(\ga)_k (\al+\be+1-\ga)_k}{(\al+1)_k (\be+1)_k}\,
(x+y)^k P_k^{(\al,\be)}\left(\frac{1+xy}{x+y}\right).
\label{199}
\end{multline}
In \myciteKLS{109}{(14)} the case $\ga$ nonpositive integer of (9.8.15)
is given. When we do this for \eqref{199} with $\ga=-n\in\ZZ_{\le0}$
this yields the
inverse of Bateman's bilinear sum, as is given in
\myciteKLS{331}{(2.19), (2.20)}, \mycite{DLMF}{(18.18.25), (18.18.26)}.
\paragraph{Bilinear generating functions}
For $0\le r<1$ and $x,y\in[-1,1]$ we have in terms of $F_4$ (see~\eqref{62}):
\begin{align}
&\sum_{n=0}^\iy\frac{(\al+\be+1)_n\,n!}{(\al+1)_n(\be+1)_n}\,r^n\,
P_n^{(\al,\be)}(x)\,P_n^{(\al,\be)}(y)
=\frac1{(1+r)^{\al+\be+1}}
\nonumber\\
&\qquad\quad\times F_4\Big(\thalf(\al+\be+1),\thalf(\al+\be+2);\al+1,\be+1;
\frac{r(1-x)(1-y)}{(1+r)^2},\frac{r(1+x)(1+y)}{(1+r)^2}\Big),
\label{58}\sLP
&\sum_{n=0}^\iy\frac{2n+\al+\be+1}{n+\al+\be+1}
\frac{(\al+\be+2)_n\,n!}{(\al+1)_n(\be+1)_n}\,r^n\,
P_n^{(\al,\be)}(x)\,P_n^{(\al,\be)}(y)
=\frac{1-r}{(1+r)^{\al+\be+2}}\nonumber\\
&\qquad\quad\times F_4\Big(\thalf(\al+\be+2),\thalf(\al+\be+3);\al+1,\be+1;
\frac{r(1-x)(1-y)}{(1+r)^2},\frac{r(1+x)(1+y)}{(1+r)^2}\Big).
\label{59}
\end{align}
Formulas \eqref{58} and \eqref{59} were first
given by Bailey \myciteKLS{91}{(2.1), (2.3)}.
See Stanton \myciteKLS{485} for a shorter proof.
(However, in the second line of
\myciteKLS{485}{(1)} $z$ and $Z$ should be interchanged.)$\;$
As observed in Bailey \myciteKLS{91}{p.10}, \eqref{59} follows
from \eqref{58}
by applying the operator $r^{-\half(\al+\be-1)}\,
\frac \dup{\dup r}\circ r^{\half(\al+\be+1)}$
to both sides of \eqref{58}.
In view of \eqref{60}, formula \eqref{59} is the Poisson kernel for Jacobi
polynomials. The \RHS\ of \eqref{59} makes clear that this kernel is positive.
See also the discussion in Askey \myciteKLS{46}{following (2.32)}.
\paragraph{Quadratic transformations}
\begin{align}
\frac{C_{2n}^{(\al+\half)}(x)}{C_{2n}^{(\al+\half)}(1)}
=\frac{P_{2n}^{(\al,\al)}(x)}{P_{2n}^{(\al,\al)}(1)}
&=\frac{P_n^{(\al,-\half)}(2x^2-1)}{P_n^{(\al,-\half)}(1)}\,,
\label{51}\\
\frac{C_{2n+1}^{(\al+\half)}(x)}{C_{2n+1}^{(\al+\half)}(1)}
=\frac{P_{2n+1}^{(\al,\al)}(x)}{P_{2n+1}^{(\al,\al)}(1)}
&=\frac{x\,P_n^{(\al,\half)}(2x^2-1)}{P_n^{(\al,\half)}(1)}\,.
\label{52}
\end{align}
See p.221, Remarks, last two formulas together with \eqref{50} and \eqref{49}.
Or see \mycite{DLMF}{(18.7.13), (18.7.14)}.
\paragraph{Differentiation formulas}
Each differentiation formula is given in two equivalent forms.
\begin{equation}
\begin{split}
\frac \dup{\dup x}\left((1-x)^\al P_n^{(\al,\be)}(x)\right)&=
-(n+\al)\,(1-x)^{\al-1} P_n^{(\al-1,\be+1)}(x),\\
\left((1-x)\frac \dup{\dup x}-\al\right)P_n^{(\al,\be)}(x)&=
-(n+\al)\,P_n^{(\al-1,\be+1)}(x).
\end{split}
\label{68}
\end{equation}
\begin{equation}
\begin{split}
\frac \dup{\dup x}\left((1+x)^\be P_n^{(\al,\be)}(x)\right)&=
(n+\be)\,(1+x)^{\be-1} P_n^{(\al+1,\be-1)}(x),\\
\left((1+x)\frac \dup{\dup x}+\be\right)P_n^{(\al,\be)}(x)&=
(n+\be)\,P_n^{(\al+1,\be-1)}(x).
\end{split}
\label{69}
\end{equation}
Formulas \eqref{68} and \eqref{69} follow from
\mycite{DLMF}{(15.5.4), (15.5.6)}
together with (9.8.1). They also follow from each other by \eqref{48}.
\paragraph{Generalized Gegenbauer polynomials}
These are defined by
\begin{equation}
S_{2m}^{(\al,\be)}(x):=\const P_m^{(\al,\be)}(2x^2-1),\qquad
S_{2m+1}^{(\al,\be)}(x):=\const x\,P_m^{(\al,\be+1)}(2x^2-1)
\label{70}
\end{equation}
in the notation of \myciteKLS{146}{p.156}
(see also \cite{K27}), while \cite[Section 1.5.2]{K26}
has $C_n^{(\la,\mu)}(x)=\const\allowbreak\times S_n^{(\la-\half,\mu-\half)}(x)$.
In \cite[Section 7.1]{K28} these polynomials are seen as special cases of the
\emph{Chihara polynomials}.
For $\al,\be>-1$ we have the orthogonality relation
\begin{equation}
\int_{-1}^1 S_m^{(\al,\be)}(x)\,S_n^{(\al,\be)}(x)\,|x|^{2\be+1}(1-x^2)^\al\,
\dup x=0\qquad(m\ne n).
\label{71}
\end{equation}
For $\be=\al-1$ generalized Gegenbauer polynomials are limit cases of
continuous $q$-ultraspherical polynomials, see \eqref{176}.

If we define the {\em Dunkl operator} $T_\mu$ by
\begin{equation}
(T_\mu f)(x):=f'(x)+\mu\,\frac{f(x)-f(-x)}x
\label{72}
\end{equation}
and if we choose the constants in \eqref{70} as
\begin{equation}
S_{2m}^{(\al,\be)}(x)=\frac{(\al+\be+1)_m}{(\be+1)_m}\, P_m^{(\al,\be)}(2x^2-1),\quad
S_{2m+1}^{(\al,\be)}(x)=\frac{(\al+\be+1)_{m+1}}{(\be+1)_{m+1}}\,
x\,P_m^{(\al,\be+1)}(2x^2-1)
\label{73}
\end{equation}
then (see \cite[(1.6)]{K5})
\begin{equation}
T_{\be+\half}S_n^{(\al,\be)}=2(\al+\be+1)\,S_{n-1}^{(\al+1,\be)}.
\label{74}
\end{equation}
Formula \eqref{74} with \eqref{73} substituted gives rise to two
differentiation formulas involving Jacobi polynomials which are equivalent to
(9.8.7) and \eqref{69}.

Composition of \eqref{74} with itself gives
\[
T_{\be+\half}^2S_n^{(\al,\be)}=4(\al+\be+1)(\al+\be+2)\,S_{n-2}^{(\al+2,\be)},
\]
which is equivalent to the composition of (9.8.7) and \eqref{69}:
\begin{equation}
\left(\frac{\dup^2}{\dup x^2}+\frac{2\be+1}x\,\frac \dup{\dup x}
\right)P_n^{(\al,\be)}(2x^2-1)
=4(n+\al+\be+1)(n+\be)\,P_{n-1}^{(\al+2,\be)}(2x^2-1).
\label{75}
\end{equation}
Formula \eqref{75} was also given in \myciteKLS{332}{(2.4)}.
\subsection*{9.8.1 Gegenbauer / Ultraspherical}
\label{sec9.8.1}
\paragraph{Notation}
Here the Gegenbauer polynomial is denoted by $C_n^\la$ instead of $C_n^{(\la)}$.
\paragraph{Orthogonality relation}
Write the \RHS\ of (9.8.20) as $h_n\,\de_{m,n}$. Then
\begin{equation}
\frac{h_n}{h_0}=
\frac\la{\la+n}\,\frac{(2\la)_n}{n!}\,,\quad
h_0=\frac{\pi^\half\,\Ga(\la+\thalf)}{\Ga(\la+1)},\quad
\frac{h_n}{h_0\,(C_n^\la(1))^2}=
\frac\la{\la+n}\,\frac{n!}{(2\la)_n}\,.
\label{61}
\end{equation}
\paragraph{Hypergeometric representation}
Beside (9.8.19) we have also
\begin{equation}
C_n^\lambda(x)=\sum_{\ell=0}^{\lfloor n/2\rfloor}\frac{(-1)^{\ell}(\lambda)_{n-\ell}}
{\ell!\;(n-2\ell)!}\,(2x)^{n-2\ell}
=(2x)^{n}\,\frac{(\lambda)_{n}}{n!}\,
\hyp21{-\thalf n,-\thalf n+\thalf}{1-\la-n}{\frac1{x^2}}.
\label{57}
\end{equation}
See \mycite{DLMF}{(18.5.10)}.
\paragraph{Special value}
\begin{equation}
C_n^{\la}(1)=\frac{(2\la)_n}{n!}\,.
\label{49}
\end{equation}
Use (9.8.19) or see \mycite{DLMF}{Table 18.6.1}.
\paragraph{Expression in terms of Jacobi}
\begin{equation}
\frac{C_n^\la(x)}{C_n^\la(1)}=
\frac{P_n^{(\la-\half,\la-\half)}(x)}{P_n^{(\la-\half,\la-\half)}(1)}\,,\qquad
C_n^\la(x)=\frac{(2\la)_n}{(\la+\thalf)_n}\,P_n^{(\la-\half,\la-\half)}(x).
\label{65}
\end{equation}
\paragraph{Re: (9.8.21)}
By iteration of recurrence relation (9.8.21):
\begin{multline}
x^2 C_n^\la(x)=
\frac{(n+1)(n+2)}{4(n+\la)(n+\la+1)}\,C_{n+2}^\la(x)+
\frac{n^2+2n\la+\la-1}{2(n+\la-1)(n+\la+1)}\,C_n^\la(x)\\
+\frac{(n+2\la-1)(n+2\la-2)}{4(n+\la)(n+\la-1)}\,C_{n-2}^\la(x).
\label{6}
\end{multline}
\paragraph{Bilinear generating functions}
\begin{multline}
\sum_{n=0}^\iy\frac{n!}{(2\la)_n}\,r^n\,C_n^\la(x)\,C_n^\la(y)
=\frac1{(1-2rxy+r^2)^\la}\,\hyp21{\thalf\la,\thalf(\la+1)}{\la+\thalf}
{\frac{4r^2(1-x^2)(1-y^2)}{(1-2rxy+r^2)^2}}\\
(r\in(-1,1),\;x,y\in[-1,1]).
\label{66}
\end{multline}
For the proof put $\be:=\al$ in \eqref{58}, then use \eqref{63} and \eqref{65}.
The Poisson kernel for Gegenbauer polynomials can be derived in a similar way
from \eqref{59}, or alternatively by applying the operator
$r^{-\la+1}\frac d{dr}\circ r^\la$ to both sides of \eqref{66}:
\begin{multline}
\sum_{n=0}^\iy\frac{\la+n}\la\,\frac{n!}{(2\la)_n}\,r^n\,C_n^\la(x)\,C_n^\la(y)
=\frac{1-r^2}{(1-2rxy+r^2)^{\la+1}}\\
\times\hyp21{\thalf(\la+1),\thalf(\la+2)}{\la+\thalf}
{\frac{4r^2(1-x^2)(1-y^2)}{(1-2rxy+r^2)^2}}\qquad
(r\in(-1,1),\;x,y\in[-1,1]).
\label{67}
\end{multline}
Formula \eqref{67} was obtained by Gasper \& Rahman \myciteKLS{234}{(4.4)}
as a limit case of their formula for the Poisson kernel for continuous
$q$-ultraspherical polynomials.
\paragraph{Trigonometric expansions}
By \mycite{DLMF}{(18.5.11), (15.8.1)}:
\begin{align}
C_n^{\la}(\cos\tha)
&=\sum_{k=0}^n\frac{(\la)_k(\la)_{n-k}}{k!\,(n-k)!}\,\eup^{\iup(n-2k)\tha}
=\eup^{\iup n\tha}\frac{(\la)_n}{n!}\,
\hyp21{-n,\la}{1-\la-n}{\eup^{-2\iup \tha}}\label{103}\\
&=\frac{(\la)_n}{2^\la n!}\,
\eup^{-\half \iup \la\pi}\eup^{\iup(n+\la)\tha}\,(\sin\tha)^{-\la}\,
\hyp21{\la,1-\la}{1-\la-n}{\frac{\iup \eup^{-\iup\tha}}{2\sin\tha}}\label{104}\\
&=\frac{(\la)_n}{n!}\,\sum_{k=0}^\iy\frac{(\la)_k(1-\la)_k}{(1-\la-n)_k k!}\,
\frac{\cos((n-k+\la)\tha+\thalf(k-\la)\pi)}{(2\sin\tha)^{k+\la}}\,.\label{105}
\end{align}
In \eqref{104} and \eqref{105} we require that
$\tfrac16\pi<\tha<\tfrac56\pi$. Then the convergence is absolute for $\la>\thalf$
and conditional for $0<\la\le\thalf$.

By \mycite{DLMF}{(14.13.1), (14.3.21), (15.8.1)]}:
\begin{align}
C_n^\la(\cos\tha)&=\frac{2\Ga(\la+\thalf)}{\pi^\half\Ga(\la+1)}\,
\frac{(2\la)_n}{(\la+1)_n}\,(\sin\tha)^{1-2\la}\,
\sum_{k=0}^\iy\frac{(1-\la)_k(n+1)_k}{(n+\la+1)_k k!}\,
\sin\big((2k+n+1)\tha\big)
\label{7}\\
&=\frac{2\Ga(\la+\thalf)}{\pi^\half\Ga(\la+1)}\,
\frac{(2\la)_n}{(\la+1)_n}\,(\sin\tha)^{1-2\la}\,
\Im\!\!\left(\eup^{\iup(n+1)\tha}\,
\hyp21{1-\la,n+1}{n+\la+1}{\eup^{2\iup\tha}}\right)\nonumber\\
&=\frac{2^\la\Ga(\la+\thalf)}{\pi^\half\Ga(\la+1)}\,
\frac{(2\la)_n}{(\la+1)_n}\,(\sin\tha)^{-\la}\,
\Re\!\!\left(\eup^{-\thalf \iup\la\pi}e^{\iup(n+\la)\tha}\,
\hyp21{\la,1-\la}{1+\la+n}{\frac{\eup^{i\tha}}{2i\sin\tha}}\right)\nonumber\\
&=\frac{2^{2\la}\Ga(\la+\thalf)}{\pi^\half\Ga(\la+1)}\,\frac{(2\la)_n}{(\la+1)_n}\,
\sum_{k=0}^\iy\frac{(\la)_k(1-\la)_k}{(1+\la+n)_k k!}\,
\frac{\cos((n+k+\la)\tha-\thalf(k+\la)\pi)}{(2\sin\tha)^{k+\la}}\,.
\label{106}
\end{align}
We require that $0<\tha<\pi$ in \eqref{7} and $\tfrac16\pi<\tha<\tfrac56\pi$ in
\eqref{106} The convergence is absolute for $\la>\thalf$ and conditional for
$0<\la\le\thalf$.
For $\la\in\Zpos$ the above series terminate after the term with
$k=\la-1$.
Formulas \eqref{7} and \eqref{106} are also given in
\mycite{Sz}{(4.9.22), (4.9.25)}.
\paragraph{Fourier transform}
\begin{equation}
\frac{\Ga(\la+1)}{\Ga(\la+\thalf)\,\Ga(\thalf)}\,
\int_{-1}^1 \frac{C_n^\la(y)}{C_n^\la(1)}\,(1-y^2)^{\la-\half}\,
\eup^{\iup xy}\,\dup y
=i^n\,2^\la\,\Ga(\la+1)\,x^{-\la}\,J_{\la+n}(x).
\label{8}
\end{equation}
See \mycite{DLMF}{(18.17.17) and (18.17.18)}.
\paragraph{Laplace transforms}
\begin{equation}
\frac2{n!\,\Ga(\la)}\,
\int_0^\iy H_n(tx)\,t^{n+2\la-1}\,\eup^{-t^2}\,\dup t=C_n^\la(x).
\label{56}
\end{equation}
See Nielsen \cite[p.48, (4) with p.47, (1) and p.28, (10)]{K4} (1918)
or Feldheim \cite[(28)]{K3} (1942).
\begin{equation}
\frac2{\Ga(\la+\thalf)}\,\int_0^1 \frac{C_n^\la(t)}{C_n^\la(1)}\,
(1-t^2)^{\la-\half}\,t^{-1}\,(x/t)^{n+2\la+1}\,\eup^{-x^2/t^2}\,\dup t
=2^{-n}\,H_n(x)\,\eup^{-x^2}\quad(\la>-\thalf).
\label{46}
\end{equation}
Use Askey \& Fitch \cite[(3.29)]{K2} for $\al=\pm\thalf$ together with
\eqref{48}, \eqref{51}, \eqref{52}, \eqref{54} and \eqref{55}.
\paragraph{Addition formula} (see \mycite{AAR}{(9.8.5$'$)]})
\begin{multline}
R_n^{(\al,\al)}\big(xy+(1-x^2)^\half(1-y^2)^\half t\big)
=\sum_{k=0}^n \frac{(-1)^k(-n)_k\,(n+2\al+1)_k}{2^{2k}((\al+1)_k)^2}\\
\times(1-x^2)^{k/2} R_{n-k}^{(\al+k,\al+k)}(x)\,(1-y^2)^{k/2} R_{n-k}^{(\al+k,\al+k)}(y)\,
\om_k^{(\al-\half,\al-\half)}\,R_k^{(\al-\half,\al-\half)}(t),
\label{108}
\end{multline}
where
\[
R_n^{(\al,\be)}(x):=P_n^{(\al,\be)}(x)/P_n^{(\al,\be)}(1),\quad
\om_n^{(\al,\be)}:=\frac{\int_{-1}^1 (1-x)^\al(1+x)^\be\,\dup x}
{\int_{-1}^1 (R_n^{(\al,\be)}(x))^2\,(1-x)^\al(1+x)^\be\,\dup x}\,.
\]
\subsection*{9.8.2 Chebyshev}
\label{sec9.8.2}
In addition to the Chebyshev polynomials $T_n$ of the first kind (9.8.35)
and $U_n$ of the second kind (9.8.36),
\begin{align}
T_n(x)&:=\frac{P_n^{(-\half,-\half)}(x)}{P_n^{(-\half,-\half)}(1)}
=\cos(n\tha),\quad x=\cos\tha,\\
U_n(x)&:=(n+1)\,\frac{P_n^{(\half,\half)}(x)}{P_n^{(\half,\half)}(1)}
=\frac{\sin((n+1)\tha)}{\sin\tha}\,,\quad x=\cos\tha,
\end{align}
we have Chebyshev polynomials $V_n$ {\em of the third kind}
and $W_n$ {\em of the fourth kind},
\begin{align}
V_n(x)&:=\frac{P_n^{(-\half,\half)}(x)}{P_n^{(-\half,\half)}(1)}
=\frac{\cos((n+\thalf)\tha)}{\cos(\thalf\tha)}\,,\quad x=\cos\tha,\\
W_n(x)&:=(2n+1)\,\frac{P_n^{(\half,-\half)}(x)}{P_n^{(\half,-\half)}(1)}
=\frac{\sin((n+\thalf)\tha)}{\sin(\thalf\tha)}\,,\quad x=\cos\tha,
\end{align}
see \cite[Section 1.2.3]{K20}. Then there is the symmetry
\begin{equation}
V_n(-x)=(-1)^n W_n(x).
\label{140}
\end{equation}

The names of Chebyshev polynomials of the third and fourth kind
and the notation $V_n(x)$ are due to Gautschi \cite{K21}.
The notation $W_n(x)$ was first used by Mason \cite{K22}.
Names and notations for Chebyshev polynomials of the third and fourth
kind are interchanged in \mycite{AAR}{Remark 2.5.3} and
\mycite{DLMF}{Table 18.3.1}.
\subsection*{9.9 Pseudo Jacobi (or Romanovski-Routh)}
\label{sec9.9}
In this section in \mycite{KLS} the pseudo Jacobi polynomial $P_n(x;\nu,N)$
in (9.9.1)
is considered
for $N\in\ZZ_{\ge0}$ and $n=0,1,\ldots,n$. However, we can more generally take
$-\thalf<N\in\RR$ (so here I overrule my convention formulated in the
beginning of this paper), $N_0$ integer such that
$N-\thalf\le N_0<N+\thalf$, and $n=0,1,\ldots,N_0$
(see \myciteKLS{382}{\S5, case A.4}). The orthogonality relation (9.9.2)
is valid for $m,n=0,1,\ldots,N_0$.
\paragraph{History}
These polynomials were first observed by Routh \cite{K13} in 1885, but not
as orthogonal polynomials (see Natanson \cite{K31} about the history).
Romanovski \myciteKLS{463} (see also Lesky \myciteKLS{382})
independently obtained them in 1929 as orthogonal
polynomials.
\paragraph{Limit relation:}
{\bf Pseudo big $q$-Jacobi $\longrightarrow$ Pseudo Jacobi}\\
See also \eqref{118}.
\paragraph{References}
See also \mycite{Ism}{\S20.1}, \myciteKLS{51},
\myciteKLS{384}, \cite{K11}, \cite{K10}, \cite{K12}.
\subsection*{9.10 Meixner}
\label{sec9.10}
\paragraph{History}
In 1934 Meixner \myciteKLS{406} (see
(1.1) and case IV on pp.~10, 11 and 12) gave the orthogonality
measure for the polynomials $P_n$ given by the generating function
\[
e^{x u(t)}\,f(t)=\sum_{n=0}^\iy P_n(x)\,\frac{t^n}{n!}\,,
\]
where
\[
e^{u(t)}=\left(\frac{1-\be t}{1-\al t}\right)^{\frac1{\al-\be}},\quad
f(t)=\frac{(1-\be t)^{\frac{k_2}{\be(\al-\be)}}}{(1-\al t)^{\frac{k_2}{\al(\al-\be)}}}\quad
(k_2<0;\;\al>\be>0\;\;{\rm or}\;\;\al<\be<0).
\]
Then $P_n$ can be expressed as a Meixner polynomial:
\[
P_n(x)=(-k_2(\al\be)^{-1})_n\,\be^n\,
M_n\left(-\,\frac{x+k_2\al^{-1}}{\al-\be},-k_2(\al\be)^{-1},\be\al^{-1}\right).
\]

In 1938 Gottlieb \cite[\S2]{K1} introduces polynomials $l_n$ ``of Laguerre type''
which turn out to be special Meixner polynomials:
$l_n(x)=\eup^{-n\la} M_n(x;1,e^{-\la})$.
\paragraph{Uniqueness of orthogonality measure}
The coefficient of $p_{n-1}(x)$ in (9.10.4) behaves as $O(n^2)$ as $n\to\iy$.
Hence \eqref{93} holds, by which the orthogonality measure is unique.
\subsection*{9.11 Krawtchouk}
\label{sec9.11}
\paragraph{Special values}
By (9.11.1) and the binomial formula:
\begin{equation}
K_n(0;p,N)=1,\qquad
K_n(N;p,N)=(1-p^{-1})^n.
\label{9}
\end{equation}
The self-duality (p.240, Remarks, first formula)
\begin{equation}
K_n(x;p,N)=K_x(n;p,N)\qquad (n,x\in \{0,1,\ldots,N\})
\label{147}
\end{equation}
combined with \eqref{9} yields:
\begin{equation}
K_N(x;p,N)=(1-p^{-1})^x\qquad(x\in\{0,1,\ldots,N\}).
\label{148}
\end{equation}
\paragraph{Symmetry}
By the orthogonality relation (9.11.2):
\begin{equation}
\frac{K_n(N-x;p,N)}{K_n(N;p,N)}=K_n(x;1-p,N).
\label{10}
\end{equation}
By \eqref{10} and \eqref{147} we have also
\begin{equation}
\frac{K_{N-n}(x;p,N)}{K_N(x;p,N)}=K_n(x;1-p,N)
\qquad(n,x\in\{0,1,\ldots,N\}),
\label{149}
\end{equation}
and, by \eqref{149}, \eqref{10} and \eqref{9},
\begin{equation}
K_{N-n}(N-x;p,N)=\left(\frac p{p-1}\right)^{n+x-N}K_n(x;p,N)
\qquad(n,x\in\{0,1,\ldots,N\}).
\label{150}
\end{equation}
A particular case of \eqref{10} is:
\begin{equation}
K_n(N-x;\thalf,N)=(-1)^n K_n(x;\thalf,N).
\label{11}
\end{equation}
Hence
\begin{equation}
K_{2m+1}(N;\thalf,2N)=0.
\label{12}
\end{equation}
From (9.11.11):
\begin{equation}
K_{2m}(N;\thalf,2N)=\frac{(\thalf)_m}{(-N+\thalf)_m}\,.
\label{13}
\end{equation}
\paragraph{Quadratic transformations}
\begin{align}
K_{2m}(x+N;\thalf,2N)&=\frac{(\thalf)_m}{(-N+\thalf)_m}\,
R_m(x^2;-\thalf,-\thalf,N),
\label{31}\\
K_{2m+1}(x+N;\thalf,2N)&=-\,\frac{(\tfrac32)_m}{N\,(-N+\thalf)_m}\,
x\,R_m(x^2-1;\thalf,\thalf,N-1),
\label{33}\\
K_{2m}(x+N+1;\thalf,2N+1)&=\frac{(\tfrac12)_m}{(-N-\thalf)_m}\,
R_m(x(x+1);-\thalf,\thalf,N),
\label{32}\\
K_{2m+1}(x+N+1;\thalf,2N+1)&=\frac{(\tfrac32)_m}{(-N-\thalf)_{m+1}}\,
(x+\thalf)\,R_m(x(x+1);\thalf,-\thalf,N),
\label{34}
\end{align}
where $R_m$ is a dual Hahn polynomial (9.6.1). For the proofs use
(9.6.2), (9.11.2), (9.6.4) and (9.11.4).
\paragraph{Recurrence relation}
Formula (9.11.3) holds for $n=N$ if we replace there the term\\
$p(N-n)K_{n+1}(x;p,N)$ by $(-x)_{N+1}/(p^N N!)$.
\paragraph{Generating functions}
\begin{multline}
\sum_{x=0}^N\binom Nx K_m(x;p,N)K_n(x;q,N)z^x\\
=\left(\frac{p-z+pz}p\right)^m
\left(\frac{q-z+qz}q\right)^n
(1+z)^{N-m-n}
K_m\left(n;-\,\frac{(p-z+pz)(q-z+qz)}z,N\right).
\label{107}
\end{multline}
This follows immediately from Rosengren \cite[(3.5)]{K8}, which goes back
to Meixner \cite{K9}.
\subsection*{9.12 Laguerre}
\label{sec9.12}
\paragraph{Notation}
Here the Laguerre polynomial is denoted by $L_n^\al$ instead of
$L_n^{(\al)}$.
\paragraph{Hypergeometric representation}
\begin{align}
L_n^\al(x)&=
\frac{(\al+1)_n}{n!}\,\hyp11{-n}{\al+1}x
\label{182}\\
&=\frac{(-x)^n}{n!} \hyp20{-n,-n-\al}-{-\,\frac1x}
\label{183}\\
&=\frac{(-x)^n}{n!}\,C_n(n+\al;x),
\label{184}
\end{align}
where $C_n$ in \eqref{184} is a
\hyperref[sec9.14]{Charlier polynomial}.
Formula \eqref{182} is (9.12.1). Then \eqref{183} follows by reversal
of summation. Finally \eqref{184} follows by \eqref{183} and \eqref{179}.
It is also the remark on top of p.244 in \mycite{KLS}, and it is essentially
\myciteKLS{416}{(2.7.10)}.
\paragraph{Uniqueness of orthogonality measure}
The coefficient of $p_{n-1}(x)$ in (9.12.4) behaves as $O(n^2)$ as $n\to\iy$.
Hence \eqref{93} holds, by which the orthogonality measure is unique.
\paragraph{Special value}
\begin{equation}
L_n^{\al}(0)=\frac{(\al+1)_n}{n!}\,.
\label{53}
\end{equation}
Use (9.12.1) or see \mycite{DLMF}{18.6.1)}.
\paragraph{Quadratic transformations}
\begin{align}
H_{2n}(x)&=(-1)^n\,2^{2n}\,n!\,L_n^{-1/2}(x^2),
\label{54}\\
H_{2n+1}(x)&=(-1)^n\,2^{2n+1}\,n!\,x\,L_n^{1/2}(x^2).
\label{55}
\end{align}
See p.244, Remarks, last two formulas.
Or see \mycite{DLMF}{(18.7.19), (18.7.20)}.
\paragraph{Fourier transform}
\begin{equation}
\frac1{\Ga(\al+1)}\,\int_0^\iy \frac{L_n^\al(y)}{L_n^\al(0)}\,
\eup^{-y}\,y^\al\,\eup^{\iup xy}\,\dup y=
\iup^n\,\frac{y^n}{(iy+1)^{n+\al+1}}\,,
\label{14}
\end{equation}
see \mycite{DLMF}{(18.17.34)}.
\paragraph{Differentiation formulas}
Each differentiation formula is given in two equivalent forms.
\begin{equation}
\frac \dup{\dup x}\left(x^\al L_n^\al(x)\right)=
(n+\al)\,x^{\al-1} L_n^{\al-1}(x),\qquad
\left(x\frac \dup{\dup x}+\al\right)L_n^\al(x)=
(n+\al)\,L_n^{\al-1}(x).
\label{76}
\end{equation}
\begin{equation}
\frac \dup{\dup x}\left(\eup^{-x} L_n^\al(x)\right)=
-\eup^{-x} L_n^{\al+1}(x),\qquad
\left(\frac \dup{\dup x}-1\right)L_n^\al(x)=
-L_n^{\al+1}(x).
\label{77}
\end{equation}
Formulas \eqref{76} and \eqref{77} follow from
\mycite{DLMF}{(13.3.18), (13.3.20)}
together with (9.12.1).
\paragraph{Generating functions}
The generating function (9.12.12)
is a limit case of the generating function \eqref{198} for Jacobi
polynomials by (9.8.16).
By (9.7.14) the generating function (9.12.12) is also a limit case of
the generating function (9.7.13) for Meixner--Pollaczek polynomials.
\paragraph{Generalized Hermite polynomials}
See \myciteKLS{146}{p.156}, \cite[Section 1.5.1]{K26}.
These are defined by
\begin{equation}
H_{2m}^\mu(x):=\const L_m^{\mu-\half}(x^2),\qquad
H_{2m+1}^\mu(x):=\const x\,L_m^{\mu+\half}(x^2).
\label{78}
\end{equation}
Then for $\mu>-\thalf$ we have orthogonality relation
\begin{equation}
\int_{-\iy}^{\iy} H_m^\mu(x)\,H_n^\mu(x)\,|x|^{2\mu}\eup^{-x^2}\,\dup x
=0\qquad(m\ne n).
\label{79}
\end{equation}
Let the Dunkl operator $T_\mu$ be defined by \eqref{72}.
If we choose the constants in \eqref{78} as
\begin{equation}
H_{2m}^\mu(x)=\frac{(-1)^m(2m)!}{(\mu+\thalf)_m}\,L_m^{\mu-\half}(x^2),\qquad
H_{2m+1}^\mu(x)=\frac{(-1)^m(2m+1)!}{(\mu+\thalf)_{m+1}}\,
 x\,L_m^{\mu+\half}(x^2)
 \label{80}
\end{equation}
then (see \cite[(1.6)]{K5})
\begin{equation}
T_\mu H_n^\mu=2n\,H_{n-1}^\mu.
\label{81}
\end{equation}
Formula \eqref{81} with \eqref{80} substituted gives rise to two
differentiation formulas involving Laguerre polynomials which are equivalent to
(9.12.6) and \eqref{76}.

Composition of \eqref{81} with itself gives
\[
T_\mu^2 H_n^\mu=4n(n-1)\,H_{n-2}^\mu,
\]
which is equivalent to the composition of (9.12.6) and \eqref{76}:
\begin{equation}
\left(\frac{\dup^2}{\dup x^2}+\frac{2\al+1}x\,\frac \dup{\dup x}
\right)L_n^\al(x^2)
=-4(n+\al)\,L_{n-1}^\al(x^2).
\label{82}
\end{equation}
\subsection*{9.13 Bessel}
\label{sec0.13}
\paragraph{Hypergeometric representation}
The constraint $n=0,1,2,\ldots,N$ can be omitted. All formulas in \S9.13
except (9.13.2) remain valid for all integer $n\ge0$. These more general
values of $n$ are even needed in the generating function (9.13.10).
\paragraph{Notation}
In the notation of Grosswald \myciteKLS{255} the \LHS\ of (9.13.1) has to be
replaced by $y_n(x;a+2)$.
\paragraph{Orthogonality relation}\quad\\
Replace the constraint \;$a<-2N-1$ \;in (9.13.2) by\;\;
$m,n=0,1,\dots,N=\lceil-(3+a)/2\rceil$.\\
Following Lesky \myciteKLS{382} the Bessel polynomials in case of
orthogonality relation (9.13.2) may be called
\emph{Romanovski--Bessel polynomials}.
\subsection*{9.14 Charlier}
\label{sec9.14}
\paragraph{Hypergeometric representation}
\begin{align}
C_n(x;a)&=\hyp20{-n,-x}-{-\,\frac1a}
\label{179}\\
&=\frac{(-x)_n}{a^n} \hyp11{-n}{x-n+1}a
\label{180}\\
&=\frac{n!}{(-a)^n}\,L_n^{x-n}(a),
\label{181}
\end{align}
where $L_n^\al(x)$ is a
\hyperref[sec9.12]{Laguerre polynomial}.
Formula \eqref{179} is (9.14.1). Then \eqref{180} follows by reversal
of the summation. Finally \eqref{181} follows by \eqref{180} and
(9.12.1). It is also the Remark on p.249 of \mycite{KLS}, and it
was earlier given in \myciteKLS{416}{(2.7.10)}.
\paragraph{Uniqueness of orthogonality measure}
The coefficient of $p_{n-1}(x)$ in (9.14.4) behaves as $O(n)$ as $n\to\iy$.
Hence \eqref{93} holds, by which the orthogonality measure is unique.
\subsection*{9.15 Hermite}
\label{sec9.15}
\paragraph{Uniqueness of orthogonality measure}
The coefficient of $p_{n-1}(x)$ in (9.15.4) behaves as $O(n)$ as $n\to\iy$.
Hence \eqref{93} holds, by which the orthogonality measure is unique.
\paragraph{Fourier transforms}
\begin{equation}
\frac1{\sqrt{2\pi}}\,\int_{-\iy}^\iy H_n(y)\,\eup^{-\half y^2}\eup^{\iup xy}
\,\dup y=
\iup^n\,H_n(x)\,\eup^{-\half x^2},
\label{15}
\end{equation}
see \mycite{AAR}{(6.1.15) and Exercise 6.11}.
\begin{equation}
\frac1{\sqrt\pi}\,\int_{-\iy}^\iy H_n(y)\,\eup^{-y^2}\,\eup^{\iup xy}\dup y=
\iup^n\,x^n\,\eup^{-\frac14 x^2},
\label{16}
\end{equation}
see \mycite{DLMF}{(18.17.35)}.
\begin{equation}
\frac{i^n}{2\sqrt\pi}\,\int_{-\iy}^\iy y^n\,\eup^{-\frac14 y^2}\,
\eup^{-\iup xy}\,\dup y=
H_n(x)\,\eup^{-x^2},
\label{17}
\end{equation}
see \mycite{AAR}{(6.1.4)}.
\subsection*{14.1 Askey--Wilson}
\label{sec14.1}
\paragraph{Symmetry}
The Askey--Wilson polynomials $p_n(x;a,b,c,d\,|\,q)$ are symmetric
in $a,b,c,d$.
\sLP
This follows from the orthogonality relation (14.1.2)
together with the value of its coefficient of $x^n$ given in (14.1.5b).
Alternatively, combine (14.1.1) with \mycite{GR}{(III.15)}.\\
As a consequence, it is sufficient to give generating function (14.1.13). Then the generating
functions (14.1.14), (14.1.15) will follow by symmetry in the parameters.
\paragraph{Basic hypergeometric representation}
In addition to (14.1.1) we have (in notation \eqref{111}):
\begin{multline}
p_n(\cos\tha;a,b,c,d\,|\, q)
=\frac{(a\eup^{-\iup\tha},b\eup^{-\iup\tha},c\eup^{-\iup\tha},
d\eup^{-\iup\tha};q)_n}
{(\eup^{-2\iup\tha};q)_n}\,\eup^{\iup n\tha}\\
\times {}_8W_7\big(q^{-n}\eup^{2\iup\tha};a\eup^{\iup\tha},b\eup^{\iup\tha},
c\eup^{\iup\tha},d\eup^{\iup\tha},q^{-n};q,q^{2-n}/(abcd)\big).
\label{113}
\end{multline}
This follows from (14.1.1) by combining (III.15) and (III.19) in
\mycite{GR}.
It is also given in \myciteKLS{513}{(4.2)}, but be aware for some slight errors.
The symmetry in $a,b,c,d$ is evident from \eqref{113}.
\paragraph{Special value and different notation}
\begin{equation}
p_n\big(\thalf(a+a^{-1});a,b,c,d\,|\, q\big)=a^{-n}\,(ab,ac,ad;q)_n\,,
\label{40}
\end{equation}
and similarly for arguments $\thalf(b+b^{-1})$, $\thalf(c+c^{-1})$ and
$\thalf(d+d^{-1})$ by symmetry of $p_n$ in $a,b,c,d$.
Formula \eqref{40} is an immediate consequence of (14.1.1).

We will also write
\begin{equation}
R_n(z;a,b,c,d\,|\,q):=\frac{p_n(\thalf(z+z^{-1});a,b,c,d\,|\,q)}
{p_n(\thalf(a+a^{-1});a,b,c,d\,|\,q)}
=\qhyp43{q^{-n},q^{n-1}abcd,az,az^{-1}}{ab,ac,ad}{q,q}.
\label{188}
\end{equation}
Here there is no longer full symmetry in $a,b,c,d$, only in $b,c,d$.
\paragraph{Trivial symmetry}
From (14.1.1) we see \myciteKLS{72}{(1.34)}
\begin{equation}
\begin{split}
p_n(x;a,b,c,d\,|\, q)&=(-1)^n p_n(-x;-a,-b,-c,-d\,|\, q),\\
R_n(z;a,b,c,d\,|\, q)&=R_n(-z;-a,-b,-c,-d\,|\, q).
\end{split}
\label{41}
\end{equation}
\paragraph{Duality}
Define parameters $\td a,\td b,\td c,\td d$ in terms of $a,b,c,d$ by
\begin{equation}
\td a=(q^{-1}abcd)^\half,\quad \td b=ab/\td a,\quad \td c=ac/\td a,\quad
\td d=ad/\td a.
\label{189}
\end{equation}
Jumping from one branch to the other branch in the square root in
the formula for $\td a$ implies that $\td a,\td b,\td c,\td d$ move to
$-\td a,-\td b,-\td c,-\td d$.
Repetition of the parameter transformation recovers
the original parameters up to a possible common multiplication of
$a,b,c,d$ by $-1$, while the branch choice for $\td a$ is irrelevant:
\begin{equation}
a=\big(q^{-1}\td a\td b\td c\td d\,\big)^\half,\quad b=\td a\td b/a,\quad
c=\td a\td c/a,\quad d=\td a\td d/a.
\label{190}
\end{equation}
From \eqref{188} we have the duality relation
\begin{equation}
R_n\big(aq^m;a,b,c,d\,|\,q\big)=
R_m\big(\td a q^n;\td a,\td b,\td c,\td d\,|\,q\big)\qquad
(m,n\in\ZZ_{\ge0}).
\label{191}
\end{equation}
By \eqref{41} both sides of \eqref{191} are invariant under
common multiplication by $-1$ of $a,b,c,d$, respectively
$\td a,\td b,\td c,\td d$.
\paragraph{Orthogonality relation}
The conditions on the parameters in (14.1.2) can be slightly relaxed:\\
Let $|a|,|b|,|c|,|d|\le1$ such that pairwise products of $a,b,c,d$ are
not equal to 1 and such that non-real parameters occur in complex
conjugate pairs.

In fact, the only possible cases which then offend the condition
$|a|,|b|,|c|,|d|<1$ are that either
precisely one parameter has absolute value 1 and equals 1 or $-1$, or
precisely two parameter values have absolute value 1, one equal to 1
and the other equal to $-1$. Then the weight fucntion will not cause
a singularity by its factors $1\pm \eup^{\iup\tha}$ and
$1\pm \eup^{-\iup\tha}$ in
the denominator, since these are compensated by the factors
$1-\eup^{2\iup\tha}$ and $1-\eup^{-2\iup\tha}$ in the numerator.

The orthogonality (14.1.3) involving discrete terms can be given
for more general parameter values as in \myciteKLS{72}{Theorem 2.5}.
There $a,b,c,d$ are real or occur in complex conjugate pairs if
non-real, and pairwise products have absolute value $\le1$ but
are not equal to 1.

\paragraph{Re: (14.1.5)}
Let
\begin{equation}
p_n(x):=\frac{p_n(x;a,b,c,d\,|\, q)}{2^n(abcdq^{n-1};q)_n}=x^n+\wt k_n x^{n-1}
+\cdots\;.
\label{18}
\end{equation}
Then
\begin{equation}
\wt k_n=-\frac{(1-q^n)(a+b+c+d-(abc+abd+acd+bcd)q^{n-1})}
{2(1-q)(1-abcdq^{2n-2})}\,.
\label{19}
\end{equation}
This follows because $\tilde k_n-\tilde k_{n+1}$ equals the coefficient
$\thalf\bigl(a+a^{-1}-(A_n+C_n)\bigr)$ of $p_n(x)$ in (14.1.5).
\paragraph{$q$-Difference equation}
The $q$-difference operator acting on $P_n(z)$ on the \RHS\ of (14.1.7),
gives, when acting on
$Q_n(z):=(az,az^{-1};q)_\iy$,
the result
\begin{multline}
q^{-n}(1-q^n)(1-abcdq^{n-1})Q_n(z)-
q^{-n}(1-abq^{n-1})(1-acq^{n-1})(1-adq^{n-1})(1-q^n)Q_{n-1}(z)\\
=A(z)Q_n(qz)-\big(A(z)+A(z^{-1})\big)Q_n(z)+A(z^{-1})Q_n(q^{-1}z).
\label{197}
\end{multline}
This formula is implicit in \cite{K32}. Use there (3.1) with the
Askey--Wilson parameters (7.15) and (7.8), and combine it with (14.1.7).
\paragraph{Generating functions}
Rahman \myciteKLS{449}{(4.1), (4.9)} gives:
\begin{align}
&\sum_{n=0}^\iy \frac{(abcdq^{-1};q)_n a^n}{(ab,ac,ad,q;q)_n}\,t^n\,
p_n(\cos\tha;a,b,c,d\,|\,q)
\nonumber\\
&=\frac{(abcdtq^{-1};q)_\iy}{(t;q)_\iy}\,
\qhyp65{(abcdq^{-1})^\half,-(abcdq^{-1})^\half,(abcd)^\half,
-(abcd)^\half,a e^{i\tha},a e^{-i\tha}}
{ab,ac,ad,abcdtq^{-1},qt^{-1}}{q,q}
\nonumber\\
&+\frac{(abcdq^{-1},abt,act,adt,ae^{i\tha},ae^{-i\tha};q)_\iy}
{(ab,ac,ad,t^{-1},ate^{i\tha},ate^{-i\tha};q)_\iy}
\nonumber\\
&\times\qhyp65{t(abcdq^{-1})^\half,-t(abcdq^{-1})^\half,t(abcd)^\half,
-t(abcd)^\half,at e^{i\tha},at e^{-i\tha}}
{abt,act,adt,abcdt^2q^{-1},qt}{q,q}\quad(|t|<1).
\label{185}
\end{align}
In the limit \eqref{109} the first term on the \RHS\ of \eqref{185}
tends to the \LHS\ of (9.1.15), while the second term tends formally
to 0. The special case $ad=bc$ of \eqref{185} was earlier given in
\myciteKLS{236}{(4.1), (4.6)}.
\paragraph{Limit relations}\quad\sLP
{\bf Askey--Wilson $\longrightarrow$ Wilson}\\
Instead of (14.1.21) we can keep a polynomial of degree $n$ while the limit is approached:
\begin{equation}
\lim_{q\to1}\frac{p_n(1-\thalf x(1-q)^2;q^a,q^b,q^c,q^d\,|\, q)}{(1-q)^{3n}}
=W_n(x;a,b,c,d).
\label{109}
\end{equation}
For the proof first derive the corresponding limit for the monic polynomials by comparing
(14.1.5) with (9.4.4).
\bLP
{\bf Askey--Wilson $\longrightarrow$ Continuous Hahn}\\
Instead of (14.4.15) we can keep a polynomial of degree $n$ while the limit is approached:
\begin{multline}
\lim_{q\uparrow1}
\frac{p_n\big(\cos\phi-x(1-q)\sin\phi;q^a \eup^{\iup\phi},
q^b \eup^{\iup\phi},q^{\overline a} \eup^{-\iup\phi},
q^{\overline b} \eup^{-\iup\phi}\,|\, q\big)}
{(1-q)^{2n}}\\
=(-2\sin\phi)^n\,n!\,p_n(x;a,b,\overline a,\overline b)\qquad
(0<\phi<\pi).
\label{177}
\end{multline}
Here the \RHS\ has a continuous Hahn polynomial (9.4.1).
For the proof first derive the corresponding limit for the monic polynomials by comparing
(14.1.5) with (9.1.5).
In fact, define the monic polynomial
\[
\wt p_n(x):=
\frac{p_n\big(\cos\phi-x(1-q)\sin\phi;q^a \eup^{\iup\phi},
q^b \eup^{\iup\phi},q^{\overline a} \eup^{-\iup\phi},
q^{\overline b} \eup^{-\iup\phi}\,|\, q\big)}
{(-2(1-q)\sin\phi)^n\,(abcdq^{n-1};q)_n}\,.
\]
Then it follows from (14.1.5) that
\begin{equation*}
x\,\wt p_n(x)=\wt p_{n+1}(x)
+\frac{(1-q^a)e^{i\phi}+(1-q^{-a})e^{-i\phi}+\wt A_n+\wt C_n}{2(1-q)\sin\phi}
\,\wt p_n(x)
+\frac{\wt A_{n-1} \wt C_n}{(1-q)^2 \sin^2\phi}\,\wt p_{n-1}(x),
\end{equation*}
where $\wt A_n$ and $\wt C_n$ are as given after (14.1.3) with
$a,b,c,d$ replaced by
$q^a \eup^{\iup\phi},q^b \eup^{\iup\phi},q^{\overline a}
\eup^{-\iup\phi},q^{\overline b} \eup^{-\iup\phi}$.
Then the recurrence equation for $\wt p_n(x)$ tends for $q\uparrow 1$ to
the recurrence equation (9.4.4) with $c=\overline a$, $d=\overline b$.
\bLP
{\bf Askey--Wilson $\longrightarrow$ Meixner--Pollaczek}\\
Instead of (14.9.15) we can keep a polynomial of degree $n$ while the limit is approached:
\begin{equation}
\lim_{q\uparrow1}
\frac{p_n\big(\cos\phi-x(1-q)\sin\phi;
q^\la \eup^{\iup\phi},0,q^\la \eup^{-\iup\phi},0\,|\, q\big)}{(1-q)^n}
=n!\,P_n^{(\la)}(x;\pi-\phi)\quad
(0<\phi<\pi).
\label{178}
\end{equation}
Here the \RHS\ has a Meixner--Pollaczek polynomial (9.7.1).
For the proof first derive the corresponding limit for the monic polynomials
by comparing (14.1.5) with (9.7.4).
In fact, define the monic polynomial
\[
\wt p_n(x):=
\frac{p_n\big(\cos\phi-x(1-q)\sin\phi;
q^\la \eup^{\iup\phi},0,q^\la \eup^{-\iup\phi},0\,|\, q\big)}
{(-2(1-q)\sin\phi)^n}\,.
\]
Then it follows from (14.1.5) that
\begin{equation*}
x\,\wt p_n(x)=\wt p_{n+1}(x)
+\frac{(1-q^\la)\eup^{\iup\phi}+(1-q^{-\la})\eup^{-\iup\phi}+\wt A_n+\wt C_n}
{2(1-q)\sin\phi}\,\wt p_n(x)
+\frac{\wt A_{n-1} \wt C_n}{(1-q)^2 \sin^2\phi}\,\wt p_{n-1}(x),
\end{equation*}
where $\wt A_n$ and $\wt C_n$ are as given after (14.1.3)
with $a,b,c,d$ replaced by
$q^\la \eup^{\iup\phi},0,q^\la \eup^{-\iup\phi},0$.
Then the recurrence equation for $\wt p_n(x)$ tends for $q\uparrow 1$ to
the recurrence equation (9.7.4).
\paragraph{References}
See also Koornwinder \cite{K7}.
\subsection*{14.2 $q$-Racah}
\label{sec14.2}
\paragraph{Symmetry}
\begin{equation}
R_n(x;\al,\be,q^{-N-1},\de\,|\, q)
=\frac{(\be q,\al\de^{-1}q;q)_n}{(\al q,\be\de q;q)_n}\,\de^n\,
R_n(\de^{-1}x;\be,\al,q^{-N-1},\de^{-1}\,|\, q).
\label{84}
\end{equation}
This follows from (14.2.1) combined with \mycite{GR}{(III.15)}.
\sLP
In particular,
\begin{equation}
R_n(x;\al,\be,q^{-N-1},-1\,|\, q)
=\frac{(\be q,-\al q;q)_n}{(\al q,-\be q;q)_n}\,(-1)^n\,
R_n(-x;\be,\al,q^{-N-1},-1\,|\, q),
\label{85}
\end{equation}
and
\begin{equation}
R_n(x;\al,\al,q^{-N-1},-1\,|\, q)
=(-1)^n\,R_n(-x;\al,\al,q^{-N-1},-1\,|\, q),
\label{86}
\end{equation}

\paragraph{Trivial symmetry}
Clearly from (14.2.1):
\begin{equation}
R_n(x;\al,\be,\ga,\de\,|\, q)=R_n(x;\be\de,\al\de^{-1},\ga,\de\,|\, q)
=R_n(x;\ga,\al\be\ga^{-1},\al,\ga\de\al^{-1}\,|\, q).
\label{83}
\end{equation}
For $\al=q^{-N-1}$ this shows that the three cases
$\al q=q^{-N}$ or $\be\de q=q^{-N}$ or $\ga q=q^{-N}$ of (14.2.1)
are not essentially different.
\paragraph{Duality}
It follows from (14.2.1) that
\begin{equation}
R_n(q^{-y}+\ga\de q^{y+1};q^{-N-1},\be,\ga,\de\,|\, q)
=R_y(q^{-n}+\be q^{n-N};\ga,\de,q^{-N-1},\be\,|\, q)\quad
(n,y=0,1,\ldots,N).
\end{equation}
\subsection*{14.3 Continuous dual $q$-Hahn}
\label{sec14.3}
The continuous dual $q$-Hahn polynomials are the special case $d=0$ of the
Askey--Wilson polynomials:
\[
p_n(x;a,b,c\,|\, q):=p_n(x;a,b,c,0\,|\, q).
\]
Hence all formulas in \S14.3 are specializations for $d=0$ of formulas in \S14.1.
\subsection*{14.4 Continuous $q$-Hahn}
\label{sec14.4}
The continuous $q$-Hahn polynomials are the special case
of Askey--Wilson polynomials with parameters
$a \eup^{\iup\phi},b \eup^{\iup\phi},a \eup^{-\iup\phi},b \eup^{-\iup\phi}$:
\[
p_n(x;a,b,\phi\,|\, q):=
p_n(x;a \eup^{\iup\phi},b \eup^{\iup\phi},a \eup^{-\iup\phi},
b \eup^{-\iup\phi}\,|\, q).
\]
In \myciteKLS{72}{(4.39)} and \mycite{GR}{(7.5.43)}
(who write $p_n(x;a,b\,|\,q)$, $x=\cos(\tha+\phi)$)
and in \mycite{KLS}{\S14.4} (who writes $p_n(x;a,b,c,d;q)$,
$x=\cos(\tha+\phi)$)
the parameter
dependence on $\phi$ is incorrectly omitted.

Since all formulas in \S14.4 are specializations of formulas in \S14.1,
there is no real need to give these specializations explicitly.
In particular, the limit (14.4.15) is in fact a limit from Askey--Wilson to
continuous Hahn. See also \eqref{177}.
\subsection*{14.5 Big $q$-Jacobi}
\label{sec14.5}
\paragraph{Different notation}
See p.442, Remarks:
\begin{equation}
P_n(x;a,b,c,d;q):=P_n(qac^{-1}x;a,b,-ac^{-1}d;q)
=\qhyp32{q^{-n},q^{n+1}ab,qac^{-1}x}{qa,-qac^{-1}d}{q,q}.
\label{123}
\end{equation}
Furthermore,
\begin{equation}
P_n(x;a,b,c,d;q)=P_n(\la x;a,b,\la c,\la d;q),
\label{141}
\end{equation}
\begin{equation}
P_n(x;a,b,c;q)=P_n(-q^{-1}c^{-1}x;a,b,-ac^{-1},1;q)
\label{142}
\end{equation}
\paragraph{Orthogonality relation}
(equivalent to (14.5.2), see also \cite[(2.42), (2.41), (2.36), (2.35)]{K17}).
Let $c,d>0$ and either $a\in (-c/(qd),1/q)$, $b\in(-d/(cq),1/q)$ or
$a/c=-\overline b/d\notin\RR$. Then
\begin{equation}
\int_{-d}^c P_m(x;a,b,c,d;q) P_n(x;a,b,c,d;q)\,
\frac{(qx/c,-qx/d;q)_\iy}{(qax/c,-qbx/d;q)_\iy}\,\dup_qx=h_n\,\de_{m,n}\,,
\label{124}
\end{equation}
where
\begin{equation}
\frac{h_n}{h_0}=q^{\half n(n-1)}\left(\frac{q^2a^2d}c\right)^n\,
\frac{1-qab}{1-q^{2n+1}ab}\,
\frac{(q,qb,-qbc/d;q)_n}{(qa,qab,-qad/c;q)_n}
\label{125}
\end{equation}
and
\begin{equation}
h_0=(1-q)c\,\frac{(q,-d/c,-qc/d,q^2ab;q)_\iy}
{(qa,qb,-qbc/d,-qad/c;q)_\iy}\,.
\label{126}
\end{equation}
\paragraph{Other hypergeometric representation and asymptotics}
\begin{align}
&P_n(x;a,b,c,d;q)
=\frac{(-qbd^{-1}x;q)_n}{(-q^{-n}a^{-1}cd^{-1};q)_n}\,
\qhyp32{q^{-n},q^{-n}b^{-1},cx^{-1}}{qa,-q^{-n}b^{-1}dx^{-1}}{q,q}
\label{138}\\
&\qquad=(qac^{-1}x)^n\,\frac{(qb,cx^{-1};q)_n}{(qa,-qac^{-1}d;q)_n}\,
\qhyp32{q^{-n},q^{-n}a^{-1},-qbd^{-1}x}{qb,q^{1-n}c^{-1}x}
{q,-q^{n+1}ac^{-1}d}
\label{132}\\
&\qquad=(qac^{-1}x)^n\,\frac{(qb,q;q)_n}{(-qac^{-1}d;q)_n}\,
\sum_{k=0}^n\frac{(cx^{-1};q)_{n-k}}{(q,qa;q)_{n-k}}\,
\frac{(-qbd^{-1}x;q)_k}{(qb,q;q)_k}\,(-1)^k q^{\half k(k-1)}(-dx^{-1})^k.
\label{133}
\end{align}
Formula \eqref{138} follows from \eqref{123} by
\mycite{GR}{(III.11)} and next \eqref{132} follows by series inversion
\mycite{GR}{Exercise 1.4(ii)}.
Formulas \eqref{138} and \eqref{133} are also given in
\mycite{Ism}{(18.4.28), (18.4.29)}.
It follows from \eqref{132} or \eqref{133} that
(see \myciteKLS{298}{(1.17)} or \mycite{Ism}{(18.4.31)})
\begin{equation}
\lim_{n\to\iy}(qac^{-1}x)^{-n} P_n(x;a,b,c,d;q)
=\frac{(cx^{-1},-dx^{-1};q)_\iy}{(-qac^{-1}d,qa;q)_\iy}\,,
\label{134}
\end{equation}
uniformly for $x$ in compact subsets of $\CC\backslash\{0\}$.
(Exclusion of the spectral points $x=cq^m,dq^m$ ($m=0,1,2,\ldots$),
as was done in \myciteKLS{298} and \mycite{Ism}, is not necessary. However,
while \eqref{134} yields 0 at these points, a more refined asymptotics
at these points is given in \myciteKLS{298} and \mycite{Ism}.)$\;$
For the proof of \eqref{134} use that
\begin{equation}
\lim_{n\to\iy}(qac^{-1}x)^{-n} P_n(x;a,b,c,d;q)
=\frac{(qb,cx^{-1};q)_n}{(qa,-qac^{-1}d;q)_n}\,
\qhyp11{-qbd^{-1}x}{qb}{q,-dx^{-1}},
\label{135}
\end{equation}
which can be evaluated by \mycite{GR}{(II.5)}.
Formula \eqref{135} follows formally from \eqref{132}, and it follows rigorously, by
dominated convergence, from \eqref{133}.
\paragraph{Symmetry}
(see \cite[\S2.5]{K17} and combine with \eqref{123}).
\begin{equation}
\frac{P_n(x;a,b,c,d;q)}{P_n(-d/(qb);a,b,c,d;q)}
=P_n(-x;b,a,d,c;q)=P_n(x;-bcd^{-1},-ac^{-1}d,c,d;q).
\end{equation}
In particular (\emph{symmetric big $q$-Jacobi polynomials}),
\begin{equation}
P_n(-x;a,a,1,1;q)=(-1)^n P_n(x;a,a,1,1;q).
\label{216}
\end{equation}
\paragraph{Special values}
\begin{align}
P_n(c/(qa);a,b,c,d;q)&=1,\\
P_n(-d/(qb);a,b,c,d;q)&=\left(-\,\frac{ad}{bc}\right)^n\,
\frac{(qb,-qbc/d;q)_n}{(qa,-qad/c;q)_n}\,,\\
P_n(c;a,b,c,d;q)&=
q^{\half n(n+1)}\left(\frac{ad}c\right)^n
\frac{(-qbc/d;q)_n}{(-qad/c;q)_n}\,,\\
P_n(-d;a,b,c,d;q)&=q^{\half n(n+1)} (-a)^n\,\frac{(qb;q)_n}{(qa;q)_n}\,.
\end{align}
\paragraph{Recurrence relation}
See (14.5.3). For $n=1,2,\ldots$:
\begin{multline}
qac^{-1}x P_n(x;a,b,c,d;q)=A_n P_{n+1}(x;a,b,c,d;q)\\
+(1-A_n-C_n)P_n(x;a,b,c,d;q)
+C_n P_{n-1}(x;a,b,c,d;q),
\label{215}
\end{multline}
where
\begin{align*}
A_n&=\frac{(1-q^{n+1}a)(1-q^{n+1}ab)(1+q^{n+1}ac^{-1}d)}
{(1-q^{2n+1}ab)(1-q^{2n+2}ab)}\,,\\
C_n&=q^{n+1}a^2c^{-1}d\,
\frac{(1-q^n)(1+q^n bcd^{-1})(1-q^nb)}
{(1-q^{2n}ab)(1-q^{2n+1}ab)}\,.
\end{align*}
For $n=0$:
\begin{multline}
qac^{-1}x P_0(x;a,b,c,d;q)=\frac{(1-qa)(1+qac^{-1}d)}{1-q^2ab}\,
P_1(x;a,b,c,d;q)\\+\frac{qa(c-d-q(bc-ad))}{c(1-q^2ab)}\,P_0(x;a,b,c,d;q).
\label{220}
\end{multline}
In \eqref{215} we have $1-A_n-C_n=0$ for $n=1,2,\ldots$
if $a=b$, $c=d$ or $ab=1$, $acd^{-1}=1$. In \eqref{220} the last term
on the right vanishes if $a=b$, $c=d$, but not if $ab=1$, $acd^{-1}=1$,
$a\ne1$.

So for symmetric big $q$-Jacobi polynomials we have
\begin{multline}
qax P_n(x;a,a,1,1;q)=
\frac{1-q^{n+1}a^2}{1-q^{2n+1}a^2}\,P_{n+1}(x;a,a,1,1;q)\\
+q^{n+1}a^2\,
\frac{1-q^n}{1-q^{2n+1}a^2}\,P_{n-1}(x;a,a,1,1;q).
\label{217}
\end{multline}
Equivalently,
\begin{equation}
x p_n(x)=
\frac{1-q^{n+1}a^2}{1-q^{2n+1}a^2}\,p_{n+1}(x)
+\frac{q^{n-1}(1-q^n)}{1-q^{2n+1}a^2}\,p_{n-1}(x),
\label{218}
\end{equation}
where $p_n(x)=(qa)^{-n}P_n(x;a,a,1,1;q)$.
\paragraph{Second order $q$-difference equation}
(see (14.5.5). Let $P_n(x)=P_n(x;a,b,c,d;q)$.
\begin{multline}
(q^{-n}-1)(1-q^{n+1}ab)P_n(x)
=qabx^{-2}(x-q^{-1}a^{-1}c)(x+q^{-1}b^{-1}d)(P_n(qx)-P_n(x))\\
+x^{-2}(x-c)(x+d)(P_n(q^{-1}x)-P_n(x)).
\label{219}
\end{multline}
\paragraph{Quadratic transformations}
(see \cite[(2.48), (2.49)]{K17} and \eqref{128}).\\
These express big $q$-Jacobi polynomials $P_m(x;a,a,1,1;q)$ in terms of little
$q$-Jacobi polynomials (see \S14.12).
\begin{align}
P_{2n}(x;a,a,1,1;q)&=\frac{p_n(x^2;q^{-1},a^2;q^2)}{p_n((qa)^{-2};q^{-1},a^2;q^2)}\,,
\label{130}\\
P_{2n+1}(x;a,a,1,1;q)&=
\frac{qax\,p_n(x^2;q,a^2;q^2)}{p_n((qa)^{-2};q,a^2;q^2)}\,.
\label{131}
\end{align}
Hence, by (14.12.1), \mycite{GR}{Exercise 1.4(ii)} and \eqref{128},
\begin{align}
P_n(x;a,a,1,1;q)&=\frac{(qa^2;q^2)_n}{(qa^2;q)_n}\,(qax)^n\,
\qhyp21{q^{-n},q^{-n+1}}{q^{-2n+1}a^{-2}}{q^2,(ax)^{-2}}
\label{136}\\
&=\frac{(q;q)_n}{(qa^2;q)_n}\,(qa)^n\,
\sum_{k=0}^{[\half n]}(-1)^k q^{k(k-1)}
\frac{(qa^2;q^2)_{n-k}}{(q^2;q^2)_k\,(q;q)_{n-2k}}\,x^{n-2k}.
\label{137}
\end{align}
\paragraph{$q$-Chebyshev polynomials}
In \eqref{123}, with $c=d=1$, the cases $a=b=q^{-\half}$ and $a=b=q^\half$ can be considered
as $q$-analogues of the Chebyshev polynomials of the first and second kind, respectively
(\S9.8.2) because of the limit (14.5.17). The quadratic relations \eqref{130}, \eqref{131}
can also be specialized to these cases. The definition of the $q$-Chebyshev polynomials
may vary by normalization and by dilation of argument. They were considered in
\cite{K18}. 
By \myciteKLS{24}{p.279} and \eqref{130}, \eqref{131}, the {\em Al-Salam-Ismail polynomials}
$U_n(x;a,b)$ ($q$-dependence suppressed) in the case $a=q$ can be expressed as
$q$-Chebyshev polynomials of the second kind:
\begin{equation*}
U_n(x,q,b)=(q^{-3} b)^{\half n}\,\frac{1-q^{n+1}}{1-q}\,
P_n(b^{-\half}x;q^\half,q^\half,1,1;q).
\end{equation*}
Similarly, by \cite[(5.4), (5.1), (5.3)]{K19} and \eqref{130}, \eqref{131}, Cigler's $q$-Chebyshev
polynomials $T_n(x,s,q)$ and $U_n(x,s,q)$
can be expressed in terms of the $q$-Chebyshev cases of \eqref{123}:
\begin{align*}
T_n(x,s,q)&=(-s)^{\half n}\,P_n((-qs)^{-\half} x;q^{-\half},q^{-\half},1,1;q),\\
U_n(x,s,q)&=(-q^{-2}s)^{\half n}\,\frac{1-q^{n+1}}{1-q}\,
P_n((-qs)^{-\half} x;q^{\half},q^{\half},1,1;q).
\end{align*}
\paragraph{Limit to Discrete $q$-Hermite I}
\begin{equation}
\lim_{a\to0} a^{-n}\,P_n(x;a,a,1,1;q)=q^n\,h_n(x;q).
\label{139}
\end{equation}
Here $h_n(x;q)$ is given by (14.28.1).
For the proof of \eqref{139} use \eqref{138}.
\paragraph{Pseudo big $q$-Jacobi polynomials}
Let $a,b,c,d\in\CC$, $z_+>0$, $z_-<0$ such that
$\tfrac{(ax,bx;q)_\iy}{(cx,dx;q)_\iy}>0$ for $x\in z_- q^\ZZ\cup z_+ q^\ZZ$.
Then $(ab)/(qcd)>0$. Assume that $(ab)/(qcd)<1$.
Let $N$ be the largest nonnegative integer such that $q^{2N}>(ab)/(qcd)$.
Then
\begin{multline}
\int_{z_- q^\ZZ\cup z_+ q^\ZZ}P_m(cx;c/b,d/a,c/a;q)\,P_n(cx;c/b,d/a,c/a;q)\,
\frac{(ax,bx;q)_\iy}{(cx,dx;q)_\iy}\,\dup_qx=h_n\de_{m,n}\\
(m,n=0,1,\ldots,N),
\label{114}
\end{multline}
where
\begin{equation}
\frac{h_n}{h_0}=(-1)^n\left(\frac{c^2}{ab}\right)^n q^{\half n(n-1)} q^{2n}\,
\frac{(q,qd/a,qd/b;q)_n}{(qcd/(ab),qc/a,qc/b;q)_n}\,
\frac{1-qcd/(ab)}{1-q^{2n+1}cd/(ab)}
\label{115}
\end{equation}
and
\begin{equation}
h_0=\int_{z_- q^\ZZ\cup z_+ q^\ZZ}\frac{(ax,bx;q)_\iy}{(cx,dx;q)_\iy}\,\dup_qx
=(1-q)z_+\,
\frac{(q,a/c,a/d,b/c,b/d;q)_\iy}{(ab/(qcd);q)_\iy}\,
\frac{\tha(z_-/z_+,cdz_-z_+;q)}{\tha(cz_-,dz_-,cz_+,dz_+;q)}\,.
\label{116}
\end{equation}
See Groenevelt \& Koelink \cite[Prop.~2.2]{K14}.
Formula \eqref{116} was first given by Slater \cite[(5)]{K15} as an evaluation
of a sum of two ${}_2\psi_2$ series.
The same formula is given in Slater \myciteKLS{471}{(7.2.6)} and in
\mycite{GR}{Exercise 5.10}, but in both cases with the same slight error,
see \cite[2nd paragraph after Lemma 2.1]{K14} for correction.
The theta function is given by \eqref{117}.
Note that
\begin{equation}
P_n(cx;c/b,d/a,c/a;q)=P_n(-q^{-1}ax;c/b,d/a,-a/b,1;q).
\label{145}
\end{equation}

In \cite{K29} the weights of the pseudo big $q$-Jacobi polynomials
occur in certain measures on the space of $N$-point configurations
on the so-called extended Gelfand-Tsetlin graph.
\subsubsection*{Limit relations}
\paragraph{Pseudo big $q$-Jacobi $\longrightarrow$ Discrete Hermite II}
\begin{equation}
\lim_{a\to\iy}i^n q^{\half n(n-1)} P_n(q^{-1}a^{-1}\iup x;a,a,1,1;q)=
\wt h_n(x;q).
\label{144}
\end{equation}
For the proof use \eqref{137} and \eqref{143}.
Note that $P_n(q^{-1}a^{-1}\iup x;a,a,1,1;q)$ is obtained from the
\RHS\ of \eqref{145} by replacing $a,b,c,d$ by
$-\iup a^{-1},\iup a^{-1},\iup ,-\iup $.
\paragraph{Pseudo big $q$-Jacobi $\longrightarrow$ Pseudo Jacobi}
\begin{equation}
\lim_{q\uparrow1}
P_n(\iup q^{\half(-N-1+\iup \nu)}x;-q^{-N-1},-q^{-N-1},q^{-N+i\nu-1};q)
=\frac{P_n(x;\nu,N)}{P_n(-\iup ;\nu,N)}\,.
\label{118}
\end{equation}
Here the big $q$-Jacobi polynomial on the \LHS\ equals
$P_n(cx;c/b,d/a,c/a;q)$ with\\
$a=\iup q^{\half(N+1-\iup \nu)}$, $b=-\iup q^{\half(N+1+\iup \nu)}$,
$c=\iup q^{\half(-N-1+\iup \nu)}$, $d=-\iup q^{\half(-N-1-\iup \nu)}$.
\subsection*{14.7 Dual $q$-Hahn}
\label{sec14.7}
\paragraph{Orthogonality relation}
More generally we have (14.7.2) with positive weights in any of the following
cases:
(i) $0<\ga q<1$, $0<\de q<1$;\quad
(ii) $0<\ga q<1$, $\de<0$;\quad
(iii) $\ga<0$, $\de>q^{-N}$;\quad
(iv) $\ga>q^{-N}$, $\de>q^{-N}$;\quad
(v) $0<q\ga<1$, $\de=0$.
This also follows by inspection of the positivity of the coefficient of
$p_{n-1}(x)$ in (14.7.4).
Case (v) yields Affine $q$-Krawtchouk in view of (14.7.13).
\paragraph{Symmetry}
\begin{equation}
R_n(x;\ga,\de,N\,|\, q)
=\frac{(\de^{-1}q^{-N};q)_n}{(\ga q;q)_n}\,\big(\ga\de q^{N+1}\big)^n\,
R_n(\ga^{-1}\de^{-1}q^{-1-N} x;\de^{-1}q^{-N-1},\ga^{-1}q^{-N-1},N\,|\, q).
\label{89}
\end{equation}
This follows from (14.7.1) combined with \mycite{GR}{(III.11)}.
\subsection*{14.8 Al-Salam--Chihara}
\label{sec14.8}
\paragraph{Standardization and notation}
The definition (14.8.1) by $q$-hypergeometric representation follows the
convention of \myciteKLS{72}{p.25} that
$Q_n(x;a,b\,|\,q)=p_n(x;a,b,0,0\,|\,q)$, where $p_n(x;a,b,c,d\,|\,q)$ is
the Askey--Wilson polynomial (14.1.1). In \mycite{Ism}{(15.1.6)}
these polynomials are notated $p_n(x;a,b\,|\,q)$, equal to
$a^n/(ab;q)_n$ times $Q_n(x;a,b\,|\,q)$ as in (14.8.1).
\paragraph{Symmetry}
The Al-Salam--Chihara polynomials $Q_n(x;a,b\,|\, q)$ are symmetric in $a,b$.
\sLP
This follows from the orthogonality relation (14.8.2)
together with the value of its coefficient of $x^n$ given in (14.8.5b).
\paragraph{Orthogonality relation}
Just as in Section 14.1
the condition $|a|,|b|<1$ on the parameters in (14.8.2) can be slightly relaxed into $|a|,|b|\le 1$, $ab\ne1$.
\subsubsection*{$q^{-1}$-Al-Salam--Chihara}
\paragraph{Re: (14.8.1)}
For $x\in\Znonneg$:
\begin{align}
Q_n(\thalf(aq^{-x}+a^{-1}q^x);&a,b\,|\, q^{-1})=
(-1)^n b^n q^{-\half n(n-1)}\left((ab)^{-1};q\right)_n
\nonumber\\
&\qquad\qquad\qquad\qquad\qquad\quad
\times\qhyp31{q^{-n},q^{-x},a^{-2}q^x}{(ab)^{-1}}{q,q^nab^{-1}}
\label{20}\\
&=(-ab^{-1})^x\,q^{-\half x(x+1)}\,\frac{(qba^{-1};q)_x}{(a^{-1}b^{-1};q)_x}\,
\qhyp21{q^{-x},a^{-2}q^x}{qba^{-1}}{q,q^{n+1}}
\label{42}\\
&=(-ab^{-1})^x\,q^{-\half x(x+1)}\,\frac{(qba^{-1};q)_x}{(a^{-1}b^{-1};q)_x}\,
p_x(q^n;ba^{-1},(qab)^{-1};q).
\label{43}
\end{align}
Formula \eqref{20} follows from the first identity in (14.8.1).
Next \eqref{42} follows from \mycite{GR}{(III.8)}.
Finally \eqref{43} gives the little $q$-Jacobi polynomials (14.12.1).
See also \myciteKLS{79}{\S3} and \cite[\S3]{K36}.
\paragraph{Orthogonality}
\begin{multline}
\sum_{x=0}^\iy
\frac{(1-q^{2x}a^{-2}) (a^{-2},(ab)^{-1};q)_x}
{(1-a^{-2}) (q,bqa^{-1};q)_x}\,
(ba^{-1})^xq^{x^2}
(Q_mQ_n)(\thalf(aq^{-x}+a^{-1}q^x);a,b\,|\, q^{-1})\\
=\frac{(qa^{-2};q)_\iy}{(ba^{-1}q;q)_\iy}\,
(q,(ab)^{-1};q)_n\,(ab)^nq^{-n^2}\,\de_{m,n}.
\label{21}
\end{multline}
The constraints for having positive weights in \eqref{21} are
$(ab)^{-1}<1$, $0<qa^{-1}b<1$. Equivalently, we are in one of the following cases:
\begin{enumerate}
\item
$a,b>0$, $ab>1$, $qa^{-1}b<1$.
\item
$a,b<0$, $ab>1$, $qa^{-1}b<1$.
\item
$a=\iup a_0$, $b=\iup b_0$, $a_0,b_0>0$, $qa_0^{-1}b_0<1$.
\item
$a=-\iup a_0$, $b=-\iup b_0$, $a_0,b_0>0$, $qa_0^{-1}b_0<1$.
\end{enumerate}

Formula \eqref{21} with constraints follows from \eqref{43} together with (14.12.2)
and the completeness of
the orthogonal system of the little $q$-Jacobi polynomials,
See also \myciteKLS{79}{\S3}. An alternative proof is given in
\myciteKLS{64}. There combine (3.82) with (3.81), (3.67), (3.40).
\paragraph{Normalized recurrence relation}
\begin{equation}
xp_n(x)=p_{n+1}(x)+\thalf(a+b)q^{-n} p_n(x)+
\tfrac14(q^{-n}-1)(abq^{-n+1}-1)p_{n-1}(x),
\label{22}
\end{equation}
where
\[
Q_n(x;a,b\,|\, q^{-1})=2^n p_n(x).
\]
\paragraph{Limit to Big $q^{-1}$-Hermite}
In \eqref{43} and \eqref{21} replace $(a,b)$ by
$(\iup b^{-\half},\iup ab^{-\half})$ with $0<aq<1$ and $b>0$. Then let
$b\downarrow0$. By (14.8.17) and (14.12.14) we arrive at
big $q^{-1}$-Hermite polynomials as duals of $q$-Bessel polynomials.
\subsection*{14.9 $q$-Meixner--Pollaczek}
\label{sec14.9}
The $q$-Meixner--Pollaczek polynomials are the special case
of Askey--Wilson polynomials with parameters
$a e^{i\phi},0,a e^{-i\phi},0$:
\[
P_n(x;a,\phi\,|\, q):=\frac1{(q;q)_n}\,
p_n(x;a \eup^{\iup\phi},0,a \eup^{-\iup\phi},0\,|\, q)\quad
(x=\cos(\tha+\phi)).
\]
In \mycite{KLS}{\S14.9} the parameter dependence on $\phi$ is
incorrectly omitted.

Since all formulas in \S14.9 are specializations of formulas in \S14.1,
there is no real need to give these specializations explicitly.
See also \eqref{178}.

There is an error in \mycite{KLS}{(14.9.6), (14.9.8)}.
Read $x=\cos(\tha+\phi)$ instead of $x=\cos\tha$.
\subsection*{14.10 Continuous $q$-Jacobi}
\label{sec14.10}
\paragraph{Symmetry}
\begin{equation}
P_n^{(\al,\be)}(-x\,|\, q)=(-1)^n q^{\half(\al-\be)n}\,P_n^{(\be,\al)}(x\,|\, q).
\label{110}
\end{equation}
This follows from \eqref{41} and (14.1.19).
\subsection*{14.10.1 Continuous $q$-ultraspherical / Rogers}
\label{sec14.10.1}
\paragraph{Re: (14.10.17)}
\begin{equation}
C_n(\cos\tha;\be\,|\, q)=
\frac{(\be^2;q)_n}{(q;q)_n}\,\be^{-\half n}\,
\qhyp43{q^{-\half n},\be q^{\half n},\be^\half \eup^{\iup\tha},
\be^\half \eup^{-\iup\tha}}
{-\be,\be^\half q^{\frac14},-\be^\half q^{\frac14}}{q^\half,q^\half},
\label{23}
\end{equation}
see \mycite{GR}{(7.4.13), (7.4.14)}.
\paragraph{Special value} (see \myciteKLS{63}{(3.23)})
\begin{equation}
C_n\big(\thalf(\be^\half+\be^{-\half});\be\,|\, q\big)
=\frac{(\be^2;q)_n}{(q;q)_n}\,\be^{-\half n}.
\end{equation}
\paragraph{Re: (14.10.21)}
(another $q$-difference equation).
Let $C_n[\eup^{\iup\tha};\be\,|\, q]:=C_n(\cos\tha;\be\,|\, q)$.
\begin{equation}
\frac{1-\be z^2}{1-z^2}\,C_n[q^\half z;\be\,|\, q]+
\frac{1-\be z^{-2}}{1-z^{-2}}\,C_n[q^{-\half}z;\be\,|\, q]=
(q^{-\half n}+q^{\half n} \be)\,C_n[z;\be\,|\, q],
\label{24}
\end{equation}
see \myciteKLS{351}{(6.10)}.
\paragraph{Re: (14.10.23)}
This can also be written as
\begin{equation}
C_n[q^\half z;\be\,|\, q]-C_n[q^{-\half}z;\be\,|\, q]=
q^{-\half n}(\be-1)(z-z^{-1})C_{n-1}[z;q\be\,|\, q].
\label{25}
\end{equation}
Two other shift relations follow from the previous two equations:
\begin{align}
(\be+1)C_n[q^\half z;\be\,|\, q]&=(q^{-\half n}+q^{\half n}\be)C_n[z;\be\,|\, q]
+q^{-\half n}(\be-1)(z-\be z^{-1})C_{n-1}[z;q\be\,|\, q],
\label{26}\\
(\be+1)C_n[q^{-\half}z;\be\,|\, q]&=(q^{-\half n}+q^{\half n}\be)C_n[z;\be\,|\, q]
+q^{-\half n}(\be-1)(z^{-1}-\be z)C_{n-1}[z;q\be\,|\, q].
\label{27}
\end{align}
\paragraph{Trigonometric representation}
(see p.473, Remarks, first formula)
\begin{equation}
C_n(\cos\tha;\be\,|\, q)=\sum_{k=0}^n
\frac{(\be;q)_k (\be;q)_{n-k}}{(q;q)_k (q;q)_{n-k}}\,\eup^{\iup(n-2k)\tha}\,.
\label{173}
\end{equation}
\paragraph{Limit for $q\downarrow-1$}
(see \myciteKLS{63}{pp.~74--75}).
By \eqref{173} and \eqref{103} we obtain
\begin{align*}
\lim_{q\uparrow1} C_{2m}(x;-q^\la\,|\,-q)&=
C_m^{\half(\la+1)}(2x^2-1)+C_{m-1}^{\half(\la+1)}(2x^2-1),\\
\lim_{q\uparrow1} C_{2m+1}(x;-q^\la\,|\,-q)&=
2x\,C_m^{\half(\la+1)}(2x^2-1).
\end{align*}
By \eqref{65} and \mycite{HTF2}{10.6(36)} this can be rewritten as
\begin{align}
\lim_{q\uparrow1} C_{2m}(x;-q^\la\,|\,-q)&=
\frac{(\la)_m}{(\half\la)_m}\, P_m^{(\half\la,\half\la-1)}(2x^2-1),
\label{174}\\
\lim_{q\uparrow1} C_{2m+1}(x;-q^\la\,|\,-q)&=
2\,\frac{(\la+1)_m}{(\half\la+1)_m}\,x\,P_m^{(\half\la,\half\la)}(2x^2-1).
\label{175}
\end{align}
By \eqref{70} the limits \eqref{174}, \eqref{175} imply that
\begin{equation}
\lim_{q\uparrow1} C_n(x;-q^\la\,|\,-q)
=\const S_n^{(\half\la,\half\la-1)}(x),
\label{176}
\end{equation}
where the \RHS\ gives a one-parameter subclass of the
generalized Gegenbauer polynomial. Note that in
\cite[Section 7.1]{K28} the generalized Gegenbauer polynomials are
also observed as fitting in the $q=-1$ Askey scheme, but the limit
\eqref{176} is not observed there. Instead, the generalized Gegenbauer
polynomials are obtained in \cite[Figure 1]{K38}
as $q\to-1$ limits of little $q$-Jacobi polynomials.
\subsection*{14.11 Big $q$-Laguerre}
\label{sec14.11}
\paragraph{Symmetry}
The big $q$-Laguerre polynomials $P_n(x;a,b;q)$ are symmetric in $a,b$.
\sLP
This follows from (14.11.1).
As a consequence, it is sufficient to give generating function (14.11.11). Then the generating
function (14.1.12) will follow by symmetry in the parameters.
\subsection*{14.12 Little $q$-Jacobi}
\label{sec14.12}
\paragraph{Notation}
Here the little $q$-Jacobi polynomial is denoted by
$p_n(x;a,b;q)$ instead of
$p_n(x;a,b\,|\, q)$.
\paragraph{Basic Hypergeometric Representation}
In addition to (14.12.1) we have (see \cite[(2.46)]{K17})
\begin{equation}
p_n(x;a,b;q)=(-qb)^{-n} q^{-\half n(n-1)}\,\frac{(qb;q)_n}{(qa;q)_n}\,
\qhyp32{q^{-n},q^{n+1}ab,qbx}{qb,0}{q,q}.
\label{186}
\end{equation}
\paragraph{Special values}
(see \cite[\S2.4]{K17}).
\begin{align}
p_n(0;a,b;q)&=1,\label{127}\\
p_n(q^{-1}b^{-1};a,b;q)&=(-qb)^{-n}\,q^{-\half n(n-1)}\,\frac{(qb;q)_n}{(qa;q)_n}\,,\label{128}\\
p_n(1;a,b;q)&=(-a)^n\,q^{\half n(n+1)}\,\frac{(qb;q)_n}{(qa;q)_n}\,.\label{129}
\end{align}
\subsection*{14.14 Quantum $q$-Krawtchouk}
\label{sec14.14}
\paragraph{$q$-Hypergeometric representation}
For $n=0,1,\ldots,N$
(see (14.14.1) and use \eqref{151}):
\begin{align}
K_n^{\rm qtm}(y;p,N;q)
&=\qhyp21{q^{-n},y}{q^{-N}}{q,pq^{n+1}}
\label{152}\\
&=(pyq^{N+1};q)_n\,
\qhyp32{q^{-n},q^{-N}/y,0}{q^{-N},q^{-N-n}/(py)}{q,q}.
\label{153}
\end{align}
\paragraph{Special values}
By \eqref{152} and \mycite{GR}{(II.4)}:
\begin{equation}
K_n^{\rm qtm}(1;p,N;q)=1,\qquad
K_n^{\rm qtm}(q^{-N};p,N;q)=(pq;q)_n.
\label{163}
\end{equation}
By \eqref{153} and \eqref{163} we have the self-duality
\begin{equation}
\frac{K_n^{\rm qtm}(q^{x-N};p,N;q)}{K_n^{\rm qtm}(q^{-N};p,N;q)}
=\frac{K_x^{\rm qtm}(q^{n-N};p,N;q)}{K_x^{\rm qtm}(q^{-N};p,N;q)}\qquad
(n,x\in\{0,1,\ldots,N\}).
\label{167}
\end{equation}
By \eqref{163} and \eqref{167} we have also
\begin{equation}
K_N^{\rm qtm}(q^{-x};p,N;q)=(pq^N;q^{-1})_x\qquad(x\in\{0,1,\ldots,N\}).
\label{169}
\end{equation}
\paragraph{Limit for $q\to1$ to Krawtchouk} (see (14.14.14) and Section \hyperref[sec9.11]{9.11}):
\begin{align}
\lim_{q\to1} K_n^{\rm qtm}(1+(1-q)x;p,N;q)&=K_n(x;p^{-1},N),
\label{161}\\
\lim_{q\to1} K_n^{\rm qtm}(q^{-x};p,N;q)&=K_n(x;p^{-1},N).
\label{164}
\end{align}
\paragraph{Quantum $q^{-1}$-Krawtchouk}
By \eqref{152}, \eqref{163}, \eqref{154} and \eqref{156}
(see also p.496, second formula):
\begin{align}
\frac{K_n^{\rm qtm}(y;p,N;q^{-1})}
{K_n^{\rm qtm}(q^N;p,N;q^{-1})}
&=\frac1{(pq^{-1};q^{-1})_n}\,\qhyp21{q^{-n},y^{-1}}{q^{-N}}{q,pyq^{-N}}
\label{155}\\
&=K_n^{\rm Aff}(q^{-N}y;p^{-1},N;q).
\label{160}
\end{align}
Rewrite \eqref{160} as
\[
K_m^{\rm qtm}(1+(1-q^{-1})qx;p^{-1},N;q^{-1})
=((pq)^{-1};q^{-1})_n\,K_n^{\rm Aff}\Big(1+(1-q)q^{-N}\big(\tfrac{1-q^N}{1-q}-x\big);p,N;q\Big).
\]
In view of \eqref{161} and \eqref{162} this tends to \eqref{10} as $q\to1$.

The orthogonality relation (14.14.2) holds with positive weights for $q>1$
if $p>q^{-1}$.
\paragraph{History}
The origin of the name of the quantum $q$-Krawtchouk polynomials
is by their interpretation
as matrix elements of irreducible corepresentations of (the quantized
function algebra of) the quantum group $SU_q(2)$ considered
with respect to its quantum subgroup $U(1)$. The orthogonality
relation and dual orthogonality relation of these polynomials
are an expression of the unitarity of these corepresentations.
See for instance \myciteKLS{343}{Section 6}.
\subsection*{14.16 Affine $q$-Krawtchouk}
\label{sec14.16}
\paragraph{$q$-Hypergeometric representation}
For $n=0,1,\ldots,N$
(see (14.16.1)):
\begin{align}
K_n^{\rm Aff}(y;p,N;q)
&=\frac1{(p^{-1}q^{-1};q^{-1})_n}\,\qhyp21{q^{-n},q^{-N}y^{-1}}{q^{-N}}{q,p^{-1}y}
\label{156}\\
&=\qhyp32{q^{-n},y,0}{q^{-N},pq}{q,q}.
\label{157}
\end{align}
\paragraph{Self-duality}
By \eqref{157}:
\begin{equation}
K_n^{\rm Aff}(q^{-x};p,N;q)=K_x^{\rm Aff}(q^{-n};p,N;q)\qquad
(n,x\in\{0,1,\ldots,N\}).
\label{168}
\end{equation}
\paragraph{Special values}
By \eqref{156} and \mycite{GR}{(II.4)}:
\begin{equation}
K_n^{\rm Aff}(1;p,N;q)=1,\qquad
K_n^{\rm Aff}(q^{-N};p,N;q)=\frac1{((pq)^{-1};q^{-1})_n}\,.
\label{165}
\end{equation}
By \eqref{165} and \eqref{168} we have also
\begin{equation}
K_N^{\rm Aff}(q^{-x};p,N;q)=\frac1{((pq)^{-1};q^{-1})_x}\,.
\label{170}
\end{equation}
\paragraph{Limit for $q\to1$ to Krawtchouk} (see (14.16.14) and Section \hyperref[sec9.11]{9.11}):
\begin{align}
\lim_{q\to1} K_n^{\rm Aff}(1+(1-q)x;p,N;q)&=K_n(x;1-p,N),
\label{162}\\
\lim_{q\to1} K_n^{\rm Aff}(q^{-x};p,N;q)&=K_n(x;1-p,N).
\label{166}
\end{align}
\paragraph{A relation between quantum and affine $q$-Krawtchouk}\quad\\
By \eqref{152}, \eqref{156}, \eqref{165} and \eqref{168}
we have for $x\in\{0,1,\ldots,N\}$:
\begin{align}
K_{N-n}^{\rm qtm}(q^{-x};p^{-1}q^{-N-1},N;q)
&=\frac{K_x^{\rm Aff}(q^{-n};p,N;q)}{K_x^{\rm Aff}(q^{-N};p,N;q)}
\label{171}\\
&=\frac{K_n^{\rm Aff}(q^{-x};p,N;q)}{K_N^{\rm Aff}(q^{-x};p,N;q)}\,.
\label{172}
\end{align}
Formula \eqref{171} is given in \cite[formula after (12)]{K24}
and \cite[(59)]{K25}.
In view of \eqref{164} and \eqref{166}
formula \eqref{172} has \eqref{149} as a limit case for
$q\to 1$.
\paragraph{Affine $q^{-1}$-Krawtchouk}
By \eqref{156}, \eqref{165},
\eqref{154} and \eqref{152} (see also p.505, first formula):
\begin{align}
\frac{K_n^{\rm Aff}(y;p,N;q^{-1})}{K_n^{\rm Aff}(q^N;p,N;q^{-1})}
&=\qhyp21{q^{-n},q^{-N}y}{q^{-N}}{q,p^{-1}q^{n+1}}
\label{158}\\
&=K_n^{\rm qtm}(q^{-N}y;p^{-1},N;q).
\label{159}
\end{align}
Formula \eqref{159} is equivalent to \eqref{160}.
Just as for \eqref{160}, it tends after suitable substitutions to
\eqref{10} as $q\to1$.

The orthogonality relation (14.16.2) holds with positive weights for $q>1$
if $0<p<q^{-N}$.
\paragraph{History}
The affine $q$-Krawtchouk polynomials were considered by Delsarte \myciteKLS{161}{Theorem~11}, \cite[(16)]{K23}
in connection with certain association schemes.
He called these polynomials generalized Krawtchouk polynomials.
(Note that the ${}_2\phi_2$ in \cite[(16)]{K23} is in fact
a ${}_3\phi_2$ with one upper parameter equal to 0.)$\;$
Next Dunkl \myciteKLS{186}{Definition 2.6, Section 5.1}
reformulated this as an interpretation as spherical functions
on certain Chevalley groups. He called these polynomials
$q$-Kratchouk polynomials. The current name
{\em affine $q$-Krawtchouk polynomials} was introduced by
Stanton \myciteKLS{488}{(4.13)}. He chose this name because,
in \myciteKLS{488}{pp.~115--116} the polynomials arise in connection
with an affine action of a group $G$ on a space $X$. Here
$X$ is the set of $(v-n)\times n$ matrices over ${\rm GF}(q)$.
Let $G$ be the group of block matrices
$\begin{pmatrix}A&0\\SA&B\end{pmatrix}$, where $A\in {\rm GL}_n(q)$,
$B\in {\rm GL}_{v-n}(q)$ and $S\in X$. Then $G$ acts on $X$ by
$\begin{pmatrix}A&0\\SA&B\end{pmatrix}\cdot T=BTA^{-1}+S$.
\subsection*{14.17 Dual $q$-Krawtchouk}
\label{sec14.17}
\paragraph{Symmetry}
\begin{equation}
K_n(x;c,N\,|\, q)=c^n\,K_n(c^{-1}x;c^{-1},N\,|\, q).
\label{87}
\end{equation}
This follows from (14.17.1) combined with \mycite{GR}{(III.11)}.
\sLP
In particular,
\begin{equation}
K_n(x;-1,N\,|\, q)=(-1)^n\,K_n(-x;-1,N\,|\, q).
\label{88}
\end{equation}
\subsection*{14.20 Little $q$-Laguerre / Wall}
\label{sec14.20}
\paragraph{Notation}
Here the little $q$-Laguerre polynomial is denoted by
$p_n(x;a;q)$ instead of
$p_n(x;a\,|\, q)$.
\paragraph{Re: (14.20.11)}
The \RHS\ of this generating function converges for $|xt|<1$.
We can rewrite the \LHS\ by use of the transformation
\begin{equation*}
\qhyp21{0,0}c{q,z}=\frac1{(z;q)_\iy}\,\qhyp01-c{q,cz}.
\end{equation*}
Then we obtain:
\begin{equation}
(t;q)_\iy\,\qhyp21{0,0}{aq}{q,xt}
=\sum_{n=0}^\iy\frac{(-1)^n\,q^{\half n(n-1)}}{(q;q)_n}\,
p_n(x;a;q)\,t^n\qquad(|xt|<1).
\label{35}
\end{equation}
\subsubsection*{Expansion of $x^n$}
Divide both sides of \eqref{35} by $(t;q)_\iy$. Then coefficients of the
same power of $t$ on both sides must be equal. We obtain:
\begin{equation}
x^n=(a;q)_n\,\sum_{k=0}^n \frac{(q^{-n};q)_k}{(q;q)_k}\,q^{nk}\,p_k(x;a;q).
\label{36}
\end{equation}
\subsubsection*{Quadratic transformations}
Little $q$-Laguerre polynomials $p_n(x;a;q)$ with $a=q^{\pm\half}$ are
related to discrete $q$-Hermite I polynomials $h_n(x;q)$:
\begin{align}
p_n(x^2;q^{-1};q^2)&=
\frac{(-1)^n q^{-n(n-1)}}{(q;q^2)_n}\,h_{2n}(x;q),
\label{28}\\
xp_n(x^2;q;q^2)&=
\frac{(-1)^n q^{-n(n-1)}}{(q^3;q^2)_n}\,h_{2n+1}(x;q).
\label{29}
\end{align}
\subsection*{14.21 $q$-Laguerre}
\label{sec14.21}
\paragraph{Notation}
Here the $q$-Laguerre polynomial is denoted by $L_n^\al(x;q)$ instead of
$L_n^{(\al)}(x;q)$.
\subsubsection*{Orthogonality relation}
(14.21.2) can be rewritten with simplified \RHS:
\begin{equation}
\int_0^\iy L_m^{\al}(x;q)\,L_n^{\al}(x;q)\,\frac{x^\al}{(-x;q)_\iy}\,
\dup x=h_n\,\de_{m,n}
\qquad(\al>-1)
\label{119}
\end{equation}
with
\begin{equation}
\frac{h_n}{h_0}=\frac{(q^{\al+1};q)_n}{(q;q)_n q^n},\qquad
h_0=-\,\frac{(q^{-\al};q)_\iy}{(q;q)_\iy}\,\frac\pi{\sin(\pi\al)}\,.
\label{120}
\end{equation}
The expression for $h_0$ (which is Askey's $q$-gamma evaluation
\cite[(4.2)]{K16})
should be interpreted by continuity in $\al$ for
$\al\in\Znonneg$.
Explicitly we can write
\begin{equation}
h_n=q^{-\half\al(\al+1)}\,(q;q)_\al\,\log(q^{-1})\qquad(\al\in\Znonneg).
\label{121}
\end{equation}
\subsubsection*{Expansion of $x^n$}
\begin{equation}
x^n=q^{-\half n(n+2\al+1)}\,(q^{\al+1};q)_n\,
\sum_{k=0}^n\frac{(q^{-n};q)_k}{(q^{\al+1};q)_k}\,q^k\,L_k^\al(x;q).
\label{37}
\end{equation}
This follows from \eqref{36} by the equality given in the Remark at the end
of \S14.20. Alternatively, it can be derived in the same way as \eqref{36}
from the generating function (14.21.14).
\subsubsection*{Quadratic transformations}
$q$-Laguerre polynomials $L_n^\al(x;q)$ with $\al=\pm\half$ are
related to discrete $q$-Hermite II polynomials $\wt h_n(x;q)$:
\begin{align}
L_n^{-1/2}(x^2;q^2)&=
\frac{(-1)^n q^{2n^2-n}}{(q^2;q^2)_n}\,\wt h_{2n}(x;q),
\label{38}\\
xL_n^{1/2}(x^2;q^2)&=
\frac{(-1)^n q^{2n^2+n}}{(q^2;q^2)_n}\,\wt h_{2n+1}(x;q).
\label{39}
\end{align}
These follows from \eqref{28} and \eqref{29}, respectively, by applying
the equalities given in the Remarks at the end of \S14.20 and \S14.28.
\subsection*{14.27 Stieltjes-Wigert}
\label{sec14.27}
\subsubsection*{An alternative weight function}
The formula on top of p.547 should be corrected as
\begin{equation}
w(x)=\frac\ga{\sqrt\pi}\,x^{-\half}\exp(-\ga^2\ln^2 x),\quad x>0,\quad
{\rm with}\quad\ga^2=-\,\frac1{2\ln q}\,.
\label{94}
\end{equation}
For $w$ the weight function given in \mycite{Sz}{\S2.7} the \RHS\ of \eqref{94}
equals $\const w(q^{-\half}x)$. See also
\mycite{DLMF}{\S18.27(vi)}.
\subsection*{14.28 Discrete $q$-Hermite I}
\label{sec14.28}
\paragraph{History}
Discrete $q$ Hermite I polynomials (not yet with this name) first occurred in
Hahn \myciteKLS{261}, see there p.29, case V and the $q$-weight $\pi(x)$ given by
the second expression on line 4 of p.30. However note that on the line on p.29 dealing with
case V, one should read $k^2=q^{-n}$ instead of $k^2=-q^n$. Then, with the indicated
substitutions,  \myciteKLS{261}{(4.11), (4.12)} yield constant multiples of
$h_{2n}(q^{-1}x;q)$ and $h_{2n+1}(q^{-1}x;q)$, respectively,
 due to the quadratic transformations \eqref{28}, \eqref{29} together with  (4.20.1).
\subsection*{14.29 Discrete $q$-Hermite II}
\label{sec14.29}
\paragraph{Basic hypergeometric representation}(see (14.29.1))
\begin{equation}
\wt h_n(x;q)=x^n\,\qhyp21{q^{-n},q^{-n+1}}0{q^2,-q^2 x^{-2}}.
\label{143}
\end{equation}
\renewcommand{\refname}{Standard references}

\renewcommand{\refname}{References from Koekoek, Lesky \& Swarttouw}

\makeatletter
\renewcommand\@biblabel[1]{[K#1]}
\makeatother
\renewcommand{\refname}{Other references}

\quad\\
\begin{footnotesize}
\begin{quote}
{T. H. Koornwinder, Korteweg-de Vries Institute, University of Amsterdam,\\
P.O.\ Box 94248, 1090 GE Amsterdam, The Netherlands;

\vspace{\smallskipamount}
email: }{\tt thkmath@xs4all.nl}
\end{quote}
\end{footnotesize}
\end{document}